\documentclass[]{scrartcl}

\usepackage[USenglish]{babel}
\usepackage{csquotes}

\usepackage[a4paper,top=27mm,bottom=20mm,inner=25mm,outer=20mm]{geometry}

\usepackage[%
  backend=bibtex,bibencoding=ascii,
  style=numeric-comp,
  giveninits=true, uniquename=init, 
  natbib=true,
  url=true,
  doi=true,
  isbn=false,
  backref=false,
  maxnames=99,
  ]{biblatex}
\usepackage{amsmath}
\allowdisplaybreaks
\numberwithin{equation}{section}
\usepackage{amssymb}
\usepackage{commath}
\usepackage{mathtools}
\usepackage{bbm}
\usepackage{nicefrac}
\usepackage{subdepth}
\usepackage[normalem]{ulem}
\usepackage{siunitx}
\usepackage{algorithm}
\usepackage{algpseudocode}

\usepackage{bbold}
\usepackage{amsthm}

\usepackage[plainpages=false,pdfpagelabels,hidelinks,unicode]{hyperref}

\usepackage{cleveref}

\usepackage{thmtools}
\usepackage{etoolbox}
\makeatletter
\patchcmd{\thmt@setheadstyle}
 {\bgroup\thmt@space}
 {\thmt@space}
 {}{}
\patchcmd{\thmt@setheadstyle}
 {\egroup\fi}
 {\fi}
 {}{}
\makeatother
\declaretheoremstyle[
  bodyfont=\normalfont\itshape,
  headformat=\NAME\ \NUMBER\NOTE,
]{myplain}
\declaretheoremstyle[
  headformat=\NAME\ \NUMBER\NOTE,
]{mydefinition}
\newcommand{\envqed}{{\lower-0.3ex\hbox{$\triangleleft$}}}
\declaretheorem[style=myplain,numberwithin=section]{theorem}
\declaretheorem[style=myplain,numberlike=theorem]{lemma}

\declaretheorem[style=myplain,numberlike=theorem]{corollary}
\declaretheorem[style=mydefinition,numberlike=theorem,qed=\envqed]{definition}
\declaretheorem[style=mydefinition,numberlike=theorem,qed=\envqed]{remark}

\usepackage{color}
\usepackage{graphicx}
\usepackage[small]{caption}
\usepackage{subcaption}
\begingroup\expandafter\expandafter\expandafter\endgroup
\expandafter\ifx\csname pdfsuppresswarningpagegroup\endcsname\relax
\else
\fi

\usepackage{booktabs}
\usepackage{rotating}
\usepackage{multirow}

\usepackage{enumitem}

\usepackage{newpxtext,newpxmath}

\let\epsilon\varepsilon
\let\phi\varphi
\let\rho\varrho

\providecommand\R{}
\renewcommand{\R}{\mathbb{R}}
\providecommand\C{}
\renewcommand{\C}{\mathbb{C}}

\renewcommand{\vec}[1]{\pmb{#1}}

\newcommand{\xmin}{x_\mathrm{min}}
\newcommand{\xmax}{x_\mathrm{max}}

\newcommand{\ii}{\mathrm{i}}
\renewcommand{\del}{\partial}

\usepackage{authblk}

\newenvironment{keywords}{\par\textbf{Key words.}}{\par}
\newenvironment{AMS}{\par\textbf{AMS subject classification.}}{\par}

\newcommand{\dd}{\mathrm{d}}

\title{Stability of the Active Flux Method in the Framework of Summation-by-Parts Operators}

\author[1]{Wasilij~Barsukow
}
\affil[1]{Institut de Mathématiques de Bordeaux (IMB), CNRS UMR 5251, 351 Cours de la Libération, 33405 Talence, France, wasilij.barsukow@math.u-bordeaux.fr}

\author[2]{Christian~Klingenberg
}
\author[2]{Lisa~Lechner
}
\affil[2]{University of Würzburg, Institute of Mathematics, Emil-Fischer-Straße 40, 97074 Würzburg, Germany}

\author[3,4]{Jan~Nordstr\"om
}
\affil[3]{Department of Mathematics, Linköping University, SE-581 83 Linköping, Sweden}
\affil[4]{Department of Mathematics and Applied Mathematics, University of Johannesburg, Auckland Park 2006, Johannesburg, South Africa}

\author[5]{Sigrun~Ortleb
}
\affil[5]{University of Kassel, Institute of Mathematics, Untere K\"onigsstra\ss e 86, 34117 Kassel, Germany}

\author[6]{Hendrik~Ranocha
}
\affil[6]{Johannes Gutenberg University Mainz, Institute of Mathematics, Staudingerweg 9, 55128 Mainz, Germany}

\begin{document}

\maketitle

\begin{abstract}
\noindent
  The Active Flux method is a numerical method for conservation laws using a combination of cell averages and point values as independent degrees of freedom, based on ideas from finite volumes and finite differences. This unusual mix has been shown to work well in many situations. We expand the theoretical justifications of the Active Flux method by analyzing it from the point of view of summation-by-parts (SBP) operators, which are routinely used to analyze finite difference, finite volume, and finite element schemes. We investigate in what type of setting the Active Flux method can be formulated using classical or degenerate SBP operators, yielding a first and novel approach for showing the energy stability of the Active Flux method. We present the analysis for the one-dimensional scalar linear advection equation with periodic boundary conditions on a uniform grid.

\end{abstract}



\section{Introduction}

The Active Flux schemes are a class of methods introduced to solve systems of hyperbolic conservation laws (see \cite{eymann11, eymann2011active, eymann13}), an extension of Scheme V from \cite{vanleer77}. They combine two types of degrees of freedom: cell averages and shared point values at the cell interfaces. The Active Flux method uses a globally continuous 
approximation, which is conceptually different to other finite volume

methods.
Yet, the conservative updates of the cell averages resemble a finite volume method and are given by the fluxes through the boundary of the cell. The Active Flux method approximates these fluxes by quadratures directly using the point values. 
This is in contrast to many finite volume schemes that use Riemann solvers to define the fluxes.

Two ways how the point values can be updated in time have emerged. Initially (e.g., in \cite{vanleer77,eymann13})
, it was proposed to use a short-time (approximate) solution of the initial-value problem (IVP) of the conservation law.
Such an approach is appropriate for many equations (see \cite{eymann2011active,barsukow18activeflux} for scalar conservation laws and for hyperbolic systems of conservation laws in one spatial dimension) and the resulting method is a one-stage method.
The solution of the IVP naturally includes upwinding which is helpful for stability upon explicit integration in time.
Since it is quite challenging to find short-time (approximate) third-order accurate solutions for multi-dimensional systems of conservation laws, it was proposed in \cite{abgrall23,abgrall23a} to complement a semi-discrete version of the cell average update with an ordinary differential equation (ODE) for the point value and to integrate both in time using a standard Runge-Kutta method \cite{abgrall2023combination,abgrall2023extensions,abgrall2025semi}. 

The aim of this paper is to provide the first energy stability analysis of a semi-discrete Active Flux scheme using the framework of summation-by-parts~(SBP) operators. The Active Flux method distinguishes itself through its low dissipation and dispersion property (\cite{roe21}). Especially in multiple space dimensions it has structure preserving properties that are unrivaled by comparable schemes: it maintains vortical structures even on coarse cells, it preserves stationary states, it is asymptotic preserving without changing the scheme (\cite{barsukow24affourier}).
SBP operators are used to mimic the integration-by-parts rule on the discrete level. Energy stability of the resulting numerical schemes is then analyzed using matrix properties of the involved SBP operators. For a good grasp on this concept, we recommend the review papers
\cite{svard2014review,fernandez2014review}. By the SBP property, the continuous energy stability analysis for a given partial differential equation can be transferred to a (semi\nobreakdash-)discrete energy stability analysis of the numerical scheme.
Thus, energy stability may be ensured for a variety of PDEs when augmented with appropriate numerical boundary conditions, e.g., via simultaneous approximation terms~(SATs) first introduced in \cite{carpenter1994timestable}.

Discrete derivative operators with an SBP property have originally been considered and developed in the context of finite difference schemes for hyperbolic and parabolic problems, see, e.g., \cite{kreiss1974finite,strand1994summation,carpenter1994timestable,olsson1995_I,olsson1995_II}, with the aim to construct high-order accurate, conservative and stable numerical methods for hyperbolic and parabolic PDEs including variable coefficient equations and nonlinear hyperbolic conservation laws. 
More recent investigations also focus on SBP operators within various popular classes of numerical schemes, e.g., finite volume schemes
 \cite{nordstrom2001finite,nordstrom2003finite}, continuous finite element \cite{hicken2016multidimensional,hicken2020entropy,abgrall2020analysisI}, discontinous Galerkin~(DG) schemes \cite{gassner2013skew,carpenter2014entropy,ortleb2017akinetic,chan2018discretely}, flux reconstruction~(FR) schemes \cite{huynh2007flux,vincent2011newclass,ranocha2016summation}, as well as meshless methods \cite{hicken2024constructing}. In addition, SBP operators based on general function spaces have been constructed in \cite{glaubitz2023summation,Glaubitz2023_2,Glaubitz2024_1,Glaubitz2024_2}.

The classical SBP framework considers central difference operators, whereas the Active Flux method uses upwinding.
The concept of upwinding was introduced directly in the SBP framework in \cite{mattsson2017diagonal} as a special case of dual-pair derivative operators \cite{dovgilovich2015high}. Upwind SBP operators can be interpreted as classical central SBP operators plus artificial dissipation \cite{svard2005steady,mattsson2007high} resulting in one-sided difference stencils.
Furthermore, they arise when coupling multiple SBP operators on subdomains of the complete spatial domain either by interface SATs~\cite{ranocha2021broadclass} or by numerical fluxes~\cite{ortleb2023stability}. The discovery of central and upwind SBP properties in various high-order schemes suggests that some form of SBP property is an essential feature in any type of provably energy stable, high-order accurate numerical method.

This work investigates the SBP properties of the semi-discrete Active Flux method for linear advection with periodic boundary conditions and analyzes its energy stability. To the best of our knowledge, this is the first attempt to prove energy stability for this method in the SBP framework.
The definition of an SBP operator includes a mass matrix (or norm matrix) in addition to the discrete derivative operator. This mass matrix is used in the discrete integration-by-parts rule and in the definition of the discrete energy norm under which energy stability is studied. Given a discrete derivative operator, first, a suitable mass matrix needs to be found for the SBP property to hold. In this work, we investigate the potential existence and form of such a mass matrix for the Active Flux method in one-dimensional periodic setting.  First, we consider the central version of the standard semi-discrete Active Flux method for the linear advection equation with periodic boundary conditions and show that the resulting discrete derivative operator admits a diagonal mass matrix resulting in a periodic diagonal-norm SBP operator. Energy stability of the Active Flux method in this periodic setting then follows in straightforward manner by SBP properties. 

Furthermore, we consider the original semi-discrete Active Flux method based on upwinding. For this scheme, a non-diagonal positive semi-definite mass matrix is found such that the pair of Active Flux methods for positive and negative advection coefficients results in a dual-pair of periodic upwind SBP operators. Our analysis shows that a diagonal mass matrix leading to an upwind SBP property of the upwind Active Flux method does not exist. Even though the mass matrix is not positive definite since its kernel contains the multiples of the all-ones vector $\vec{1}$, we prove energy stability of the Active Flux method in the SBP framework by an additional study of the solution component in the kernel of the norm matrix. 

In addition, nullspace consistency of the Active Flux method is investigated in this work. Nullspace consistency was introduced in~\cite{svard2019convergence} and signifies that the nullspace of the continuous derivative operator is correctly transferred to its discrete counterpart and that no spurious modes lie in the discrete kernel. Classical central SBP operators for periodic problems are not nullspace consistent due to their skew-symmetric character and the central Active Flux scheme is no exception. However, periodic upwind SBP operators can be nullspace consistent and we prove nullspace consistency of the original semi-discrete Active Flux method based on upwinding.

This paper is organized as follows. We start by introducing the Active Flux method and defining periodic SBP operators in \cref{sec:overview}. Section~\ref{sec:analysis} analyzes the Active Flux method in the SBP framework. First, a central version of the Active Flux method is considered, which fits the more classical central SBP framework. Then, the semi-discrete upwind Active Flux method is analyzed using the insight gained from the analysis of the central scheme. The stability results are summarized in \cref{sec:stability} and verified in \cref{sec:experiments} by numerical experiments. Nullspace consistency of the upwind Active Flux difference operators is proven in \cref{sec:nullspace consistency}. In order to increase readability, \cref{sec:circulant-proof-of-lem:central-AF-banded-mass} collects some detailed technical aspects employed to study the SBP properties and \cref{sec:von-Neumann-stability} provides a relation to the von Neumann stability concept.

\section{General aspects of Active Flux and SBP~operators} \label{sec:overview}
We begin by introducing the \emph{Active Flux method} for conservation laws, originally proposed in \cite{vanleer77,eymann13}. Even though this method shows its real strength in two and three space dimensions, we restrict ourselves to one space dimension for simplicity and clarity. Consider the scalar conservation law
\begin{align}
	\del_t u(x,t) + \del_x f(u(x,t)) = 0
    \label{eq:conservation_law}
\end{align}
on domain $\Omega=[x_{\min}, x_{\max}]\subset \mathbb R$ with $u: \Omega \times [0, \infty) \rightarrow \mathbb R$, $f: \mathbb R \rightarrow \mathbb R$ and periodic boundary conditions. Later, \eqref{eq:conservation_law} will be simplified to the linear constant coefficient advection equation 
\begin{align}\label{eq:linad}
	\del_t u + a \del_x u = 0. 
\end{align}
Our results hold for general constants 
$a\in\mathbb R$ but as an indicator, we often choose $a = 1$ or $a = -1$.

\subsection{The Active Flux method}\label{sec:overview-af}
\begin{figure}[htbp]
  \centering
  \includegraphics[scale=0.15]{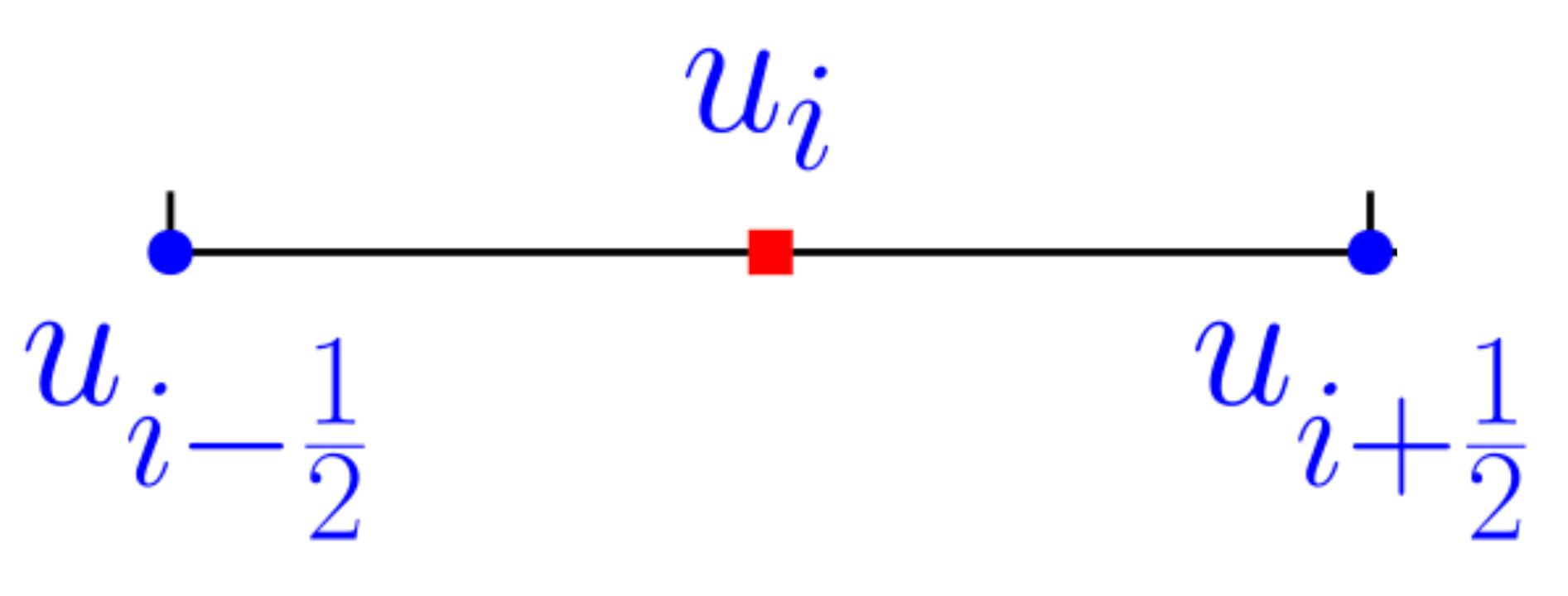}
  \caption{A cell that goes from $x_{i-\frac12}$ to $x_{i+\frac12}$ is shown: $u_i$ is the cell average of the solution $u$, whereas $u_{i+\frac12}$ and $u_{i-\frac12}$ are its point values. With these three pieces of information one can reconstruct a parabola in this cell. Since the point values are shared with neighboring cells, one associates with each cell two degrees of freedom, namely $u_{i-\frac12}$ and $u_i$. }
  \label{fig:interval}
\end{figure}

Consider a computational grid with $n$ cells $[x_{i-\frac12}, x_{i+\frac12}]$ and cell centers $x_i = \frac{x_{i+\frac12} + x_{i-\frac12}}{2}$, $i = 0, \ldots, n-1$. For simplicity, all cells have the same size $\Delta x := x_{i+\frac12} - x_{i-\frac12}$.
The Active Flux method
uses cell averages $u_i$ and point values $u_{i+\frac12}$ as degrees of freedom, see Figure~\ref{fig:interval}. The point values are placed at the cell interfaces and shared by adjacent cells.

  The three pieces of information accessible to each cell (two point values and the cell average) allow for a parabolic reconstruction $u_{\text{recon},i} \colon \left[-\frac{\Delta x}{2}, \frac{\Delta x}{2} \right ] $,
which fulfills
\begin{align}
u_{\text{recon},i}\left(\pm\frac{\Delta x}{2}\right) &= u_{i\pm\frac12}, \qquad \frac{1}{\Delta x} \int_{-\frac{\Delta x}{2}}^{\frac{\Delta x}{2}} u_{\text{recon},i}(x) \,\dd x = u_i.
\end{align}
We find $u_{\text{recon},i} = \frac{6 u_i - u_{i+\frac12} - u_{i-\frac12}}{4} + \frac{u_{i+\frac12} - u_{i-\frac12}}{\Delta x} x + 3\frac{u_{i+\frac12} + u_{i-\frac12} - 2 u_i}{\Delta x^2} x^2$.
The global reconstruction $u_\text{recon} (x) \colon \mathbb R \to \mathbb R$ is piecewise parabolic and continuous:
\begin{align}
u_\text{recon}(x) = u_{\text{recon},j}(x - x_j)
\quad \text{if} \quad x \in [x_{j-\frac12}, x_{j+\frac12}], \quad j \in \mathbb{N}.
\end{align}

When integrating the conservation law \eqref{eq:conservation_law}
over the cell, Gauss' law allows us to couple the cell averages and the point values:
\begin{align}
\frac{\dd}{\dd t} u_i + \frac{f(u_{i+\frac12}) - f(u_{i-\frac12})}{\Delta x} &= 0.\label{eq:avgupdategeneral}
\end{align}
Next, we need to determine the point value updates. In this paper, we restrict ourselves to one \emph{upwind} and one \emph{central} point update.

The semi-discrete upwind point update is given by
\begin{align}
\frac{\dd}{\dd t}u_{i+\frac12} + f'(u_{i+\frac12})(d^\text{upw}_x u)_{i+\frac12} &= 0 \label{eq:pointvalueupdategeneral}
\end{align}
combined with an upwind finite difference formula $(d^\text{upw}_x u)_{i+\frac12}$:
\begin{align}
(d^\text{upw}_x u)_{i+\frac12} := \begin{cases} \frac{\dd}{\dd x} u_\text{recon}(x_{i+\frac12}^-),  &  f'(u_{i+\frac12}) > 0,  \\\frac{\dd}{\dd x} u_\text{recon}(x_{i+\frac12}^+),   &  \text{else},   \end{cases} = \begin{cases} \frac{\dd}{\dd x} u_{\text{recon},i}\left(\frac{\Delta x}{2} \right),  & \mbox{if} \quad f'(u_{i+\frac12}) > 0,  \\\frac{\dd}{\dd x} u_{\text{recon},i+1}\left(-\frac{\Delta x}{2} \right),   &  \text{else}.   \end{cases} \label{eq:findiffgeneralupwind}
\end{align}
By construction, the finite differences are thus exact for parabolae.
The discontinuity in the derivative of $u_\text{recon}$ at $x_{i+\frac12}$ allows one to include upwinding. In \cite{abgrall23}, the scheme \eqref{eq:avgupdategeneral}-\eqref{eq:findiffgeneralupwind} has been found to be $L^2$-stable with RK3 up to a CFL number of approximately 0.4 using Fourier analysis.

From now on, we focus on the linear advection equation \eqref{eq:linad}.
For $a = 1$, we have
\begin{equation}
\label{eq:upwind-AF-Dm}
  \frac{\dif}{\dif t} u_{i+1/2} + \frac{2 u_{i-1/2} - 6 u_{i} + 4 u_{i+1/2}}{\Delta x} = 0.
\end{equation}
For $a = -1$, we have
\begin{equation}
\label{eq:upwind-AF-Dp}
   \frac{\dif}{\dif t} u_{i+1/2} - \frac{-4 u_{i+1/2} + 6 u_{i+1} - 2 u_{i+3/2}}{\Delta x} = 0.
\end{equation}
More about point updates in the Active Flux method, and convergence analyses confirming experimentally the third order of accuracy can be found in \cite{abgrall23} and papers referenced therein. The semi\nobreakdash-discrete Active Flux method can be extended to multiple dimensions and to higher orders of accuracy, see \cite{Lechner25}.

The point value update employs a \emph{central discretization} 
will be used in the theoretical analysis below. A central derivative is obtained by taking the average of the two options \eqref{eq:upwind-AF-Dm} and  \eqref{eq:upwind-AF-Dp}
\begin{align}
(d^\text{central}_x u)_{i+\frac12} := \frac12 \left( \frac{\dd}{\dd x} u_{\text{recon},i}\left(\frac{\Delta x}{2} \right)  + \frac{\dd}{\dd x} u_{\text{recon},i+1}\left(-\frac{\Delta x}{2} \right)   \right ),\label{eq:findiffgeneralcentral}
\end{align}
which leads to the following update formula
\begin{equation}
\label{eq:central_AF_point_values}
   \frac{\dif}{\dif t}  u_{i+1/2} + \frac{u_{i-1/2} - 3 u_i + 3 u_{i+1} - u_{i+3/2}}{\Delta x} = 0.
\end{equation}

\begin{remark} The Active Flux method described above is  conservative, since summing up the cell-averages \eqref{eq:avgupdategeneral} leads to conservation. 
That the point value updates inserted into this formula may themselves be found using a non-conservative update, does not change this fact.
Indeed   \eqref{eq:avgupdategeneral} will guarantee by the Lax-Wendroff theorem that if the solution converges, it converges to a weak solution.
\end{remark}

\subsection{The summation-by-parts approach for showing energy stability}\label{sec:overview-sbp}

In the following, one-dimensional periodic first-derivative SBP operators are introduced.
Consider a semi-discretization of the linear advection equation \eqref{eq:linad} with $a=1$ supplemented by periodic boundary conditions
\begin{equation}
  \frac{\dif}{\dif t} \vec{u} + D \vec{u} = \vec{0},
\end{equation}
where $\vec{u}$ is the vector of unknowns and $D$ the corresponding difference operator.
A SBP operator on the interval $\Omega = [\xmin, \xmax]$ is associated to a vector $\vec{x}$ of discrete grid nodes within $\Omega$ corresponding to the unknowns.
Due to the matrix-vector formulation of SBP schemes, we use the notation $\vec{1} = (1, \dots, 1)^T$, $\vec{0} = (0, \dots, 0)^T$ for the values of these constant functions on the given grid. \\

An SBP operator $D$ mimics the continuous integration-by-parts property on $\Omega$
\begin{align}\label{eq:ibp}
	\int_{x_{\min}}^{x_{\max}} v u_x \dd x + \int_{x_{\min}}^{x_{\max}} v_x u \dd x = vu\vert_{x_{\min}}^{x_{\max}}
\end{align}
on the discrete level. Here, we restrict ourselves to periodic boundary conditions, such that the left-hand side of \eqref{eq:ibp} vanishes.
Then, the periodic SBP property is given as
\begin{align}
	&(\mathbf v, D\mathbf u)_M + (D\mathbf v, \mathbf u)_M = 0,\label{eq:sbp}
\end{align}
where $M$ is the symmetric and positive definite norm matrix of the given SBP operator $D$ with the inner product $(u,v)_M=\mathbf u^T M \mathbf v$.

Considering \eqref{eq:linad} with $a=1$ and the periodicity of the problem we find
\begin{align}
	\frac{\dd}{\dd t}\Vert u \Vert_{L^2(\Omega)}^2 &= 2  \int_{\xmin}^{\xmax} u \partial_t u \dd x
	 = -2 \int_{\xmin}^{\xmax} u \partial_x u \dd x
= -u^2(\xmax, t)+u^2(\xmin, t) = 0\label{eq:intro_continous_energy}
\end{align}
from \eqref{eq:ibp}. Using the SBP operator $D$ with \eqref{eq:sbp}, the corresponding semi-discrete operation yields
\begin{equation}
\begin{aligned}
  \frac{\dif}{\dif t} \| \vec{u} \|_M^2
  &=
  2 \vec{u}^T M  \frac{\dif}{\dif t} \vec{u}
  =
  -2 \vec{u}^T M D \vec{u}
  =
  - \vec{u}^T M D \vec{u} + \vec{u}^T D^T M \vec{u}
  =
  0,\label{eq:intro_discrete_energy}
\end{aligned}
\end{equation}
and energy stability is ensured.
A periodic SBP operator approximating $\partial_x$ is defined as follows:
\begin{definition} \label{def:periodicSBP}
  A \emph{periodic SBP operator} on the interval
  $[\xmin, \xmax]$ consists of a grid $\vec{x}$,
  a consistent derivative matrix $D$ (with $D\vec{1}=\vec{0}$), and a symmetric and
  positive definite mass matrix~$M$ (also denoted as norm matrix) such that the periodic
  SBP property
  \begin{equation} \label{eq:SBPproperty}
    M D + D^T M = 0
  \end{equation}
holds.
Since $\int_{\xmin}^{\xmax}\dd x=\xmax - \xmin$, $M$ is scaled such that $\vec{1}^T M \vec{1} = \xmax - \xmin$.
 We call this SBP operator a \emph{diagonal-norm} SBP operator if $M$ is diagonal.
\end{definition}

\begin{remark}
As seen from Definition \ref{def:periodicSBP}, the periodic SBP property implies skew-symmetry of $D$ with respect to $M$.
\end{remark}

A classical first-derivative SBP operator $D$ results from central differences and fulfills
\eqref{eq:SBPproperty} in the periodic case, where the norm matrix $M$ assumes a role similar to quadrature rules in finite element discretizations.
Upwind SBP operators 
are a special case of dual-pair operators.
For periodic boundary conditions, dual-pair upwind SBP operators $D_+,D_-$ discretizing the first derivative fulfill a similar property: $MD_+ +(D_-)^T M = 0$. However, each of the operators $D_+$ and $D_-$ regarded separately is based on upwind differencing which leads to a built-in artificial dissipation when combined with flux splitting \cite{lundgren2020anefficient,stiernstroem2021aresidualbased,ranocha2025ontherobustness}. This approach stabilizes linear problems while retaining the stability properties associated with SBP schemes.
Here, we use the following definition:
\begin{definition}
  A \emph{periodic upwind SBP operator} on the interval
  $[\xmin, \xmax]$ consists of a grid $\vec{x}$,
  a pair of consistent derivative matrices $D_-$ and $D_+$, and a symmetric and
  positive definite mass matrix $M$ such that the periodic
  upwind SBP properties
  \begin{equation}
    M D_+ + D_-^T M = 0,
    \qquad
    M (D_+ - D_-) \text{ is negative semidefinite}
  \end{equation}
hold.
Again, $M$ is scaled
to fulfill $\vec{1}^T M \vec{1} = \xmax - \xmin$. If $M$ is diagonal, we speak of a \emph{diagonal-norm} upwind SBP operator.
\end{definition}

Using periodic upwind SBP operators, a semi-discretization
of the linear advection equation~\eqref{eq:linad} with $a=1$ is
\begin{equation}
   \frac{\dif}{\dif t}  \vec{u} + D_- \vec{u} = \vec{0}.
\end{equation}
This method is stable since
\begin{equation}
\label{eq:upwind-SBP-stability}
\begin{aligned}
  \frac{\dif}{\dif t} \| \vec{u} \|_M^2
  &=
  2 \vec{u}^T M  \frac{\dif}{\dif t} \vec{u}
  =
  -2 \vec{u}^T M D_- \vec{u}
  =
  - \vec{u}^T M D_- \vec{u} + \vec{u}^T D_+^T M \vec{u}
    =
   \vec{u}^T M (D_+ - D_-) \vec{u}
  \le
  0.
\end{aligned}
\end{equation}
Energy stability of the scheme
$\frac{\dif}{\dif t} \vec{u} - D_+ \vec{u} = \vec{0}$
for
$a = -1$ is shown similarly.

\section{Matrix analysis of the Active Flux method}\label{sec:analysis}

The following section is dedicated to the analysis of the semi-discrete Active Flux method in the framework of SBP operators.
The difficult part consists in finding a suitable mass matrix $M$ whereas the discrete derivative operator $D$ corresponding to the Active Flux method is given. For this purpose, the Active Flux scheme
\eqref{eq:avgupdategeneral} and the point update~\eqref{eq:pointvalueupdategeneral}
for the
linear scalar advection equation \eqref{eq:linad}
must be rewritten in matrix-vector form. Since periodic boundary conditions are considered,
we identify $u_{n-\frac12} := u_{-\frac12}$. In order to prove SBP properties, we will start with a central version of the Active Flux method since its SBP properties are easier to detect.
Afterwards, in the subsequent section, we will turn to the standard upwind version of the Active Flux method.

\subsection{Central Active Flux discretization}
\label{subsec:analysis_central}

As motivated before, we aim at formulating the central version of the
semi-discrete Active Flux method (given by \eqref{eq:avgupdategeneral} and \eqref{eq:findiffgeneralcentral}) within the framework of SBP operators which immediately yields stability. The first step to reach this goal consists in establishing a matrix-vector formulation.

We begin by recalling that for the cell averages, the semi-discrete Active Flux method yields the update
\begin{equation}
   \frac{\dif}{\dif t}  u_i + \frac{u_{i+1/2} - u_{i-1/2}}{\Delta x} = 0\,.
\end{equation}
We also recall that the point values are updated using a central discretization according to \eqref{eq:findiffgeneralcentral}, which yields
\begin{equation}
   \frac{\dif}{\dif t}  u_{i+1/2} + \frac{u_{i-1/2} - 3 u_i + 3 u_{i+1} - u_{i+3/2}}{\Delta x} = 0.
\end{equation}
Collecting the cell averages $u_{i}$ and the
point values $u_{i+1/2}$ in a single vector $\vec{u} \in \mathbb R^{2n}$,
we can write the semi-discrete system in matrix-vector form as
\begin{equation}
\label{eq:central-AF}
  \frac{\dif}{\dif t} \underbrace{\left(\begin{array}{c}
    \vdots \\
    u_{i-1} \\\hline
    u_{i-1/2} \\
    u_{i} \\\hline
    u_{i+1/2} \\
    u_{i+1} \\\hline
    \vdots
  \end{array}\right)}_{=:\, \vec{u}}
  + \underbrace{\frac{1}{\Delta x} \left(\begin{array}{cc|cc|cc|c}
      \ddots &&&\\
      -1 & 0 & 1 &&\\\cline{1-6}
      1 & -3 & 0 & 3 & -1 &&\\
      && -1 & 0 & 1 &&\\\hline
      && 1 & -3 & 0 & 3 & -1 \\
      &&&& -1 & 0 & 1 \\\cline{3-7}
      \multicolumn{4}{l|}{}&&& \ddots
  \end{array}\right)}_{=: \, D}
  \vec{u}
  = \vec{0}
\end{equation}
and directly obtain the discrete derivative operator $D$.
Here and below, we add guiding lines to the matrices to make it easier to see which entries belong to point values and cell averages.

We observe that $D$ is skew-symmetric with
respect to the diagonal mass matrix
\begin{equation}
\label{eq:central-AF-diagonal-mass}
  M = \frac{\Delta x}{4}
\left(\begin{array}{cc|cc|cc|c}
      \ddots &&\\
       & 3 &  \\\cline{1-4}
       &  & 1 &  &  \\
      &&  & 3 &  \\\cline{3-6}
      \multicolumn{4}{l|}{} & 1 & &  \\
      \multicolumn{4}{l|}{} && 3 &  \\\cline{5-7}
      \multicolumn{6}{l|}{}& \ddots
  \end{array}\right)
\end{equation}
since
\begin{equation}
  M D = -D^T M.
\end{equation}
Therefore, since $M$ is diagonal with positive entries, the discretization of $\partial_x$ via the Active Flux method yields a periodic SBP operator according to Definition~\ref{def:periodicSBP}.
Hence, we have proved the following theorem.

\begin{theorem}
    The central version \eqref{eq:findiffgeneralcentral}
    of the semi-discrete Active Flux method can be
    formulated using SBP operators with derivative matrix
    $D$ \eqref{eq:central-AF} and diagonal norm/mass matrix
    $M$ \eqref{eq:central-AF-diagonal-mass},
    satisfying $\vec{1}^T M \vec{1} = \xmax - \xmin$.
\end{theorem}

\begin{remark} \label{rem:trapezoidalinterpretation}
The diagonal mass matrix $M$ in \eqref{eq:central-AF-diagonal-mass} can be interpreted as a chained trapezoidal rule in each cell. Using the transformation of the cell averages $u_i = \frac 16 u_{i - \frac 12} + \frac 23 u_{i,p} + \frac 16 u_{i + \frac 12}$ to the cell midpoints $u_{i,p}$ derived from the reconstruction $u_{\text{recon},i}$ and inserting this for the terms $\frac 34 u_i$ shows
\begin{align*}
	\vec{1}^T M \vec{u} 
    &= \Delta x \left(\dots + \frac 14 u_{i - \frac 12} + \left(\frac 18 u_{i - \frac 12} + \frac 12 u_{i,p} + \frac 18 u_{i + \frac 12}\right) + \frac14u_{i+1/2} + \left(\frac 18 u_{i + \frac 12} + \frac 12 u_{i+1,p} + \frac 18 u_{i + \frac 32}\right) + \dots\right).
\end{align*}

Summing up the point values at cell interfaces and midpoints, we obtain the two variants
\begin{align*}     \vec{1}^T M \vec{u} &= \frac{\Delta x}{2} \left(\dots + u_{i - \frac 12} + u_{i,p} + u_{i + \frac 12} + u_{i+1,p} + \dots\right)\\
              &= \frac{\Delta x}{2} \left(\dots + \left(\frac 12 u_{i - \frac 12} + u_{i,p} + \frac 12 u_{i + \frac 12}\right)+ \left(\frac 12 u_{i + \frac 12} + u_{i+1,p} + \frac 12 u_{i + \frac 32}\right) \dots\right).
\end{align*}
Thus, the quadrature induced by the diagonal mass matrix $M$ defined by \eqref{eq:central-AF-diagonal-mass} can be considered either as the chained trapezoidal rule on the complete computational domain or in each individual cell. In fact, transferred to the interface and midpoint nodes $u_i, u_{i,p}$ this quadrature rule directly corresponds to the diagonal mass matrix of the classical second-order SBP operator on equidistant nodes with grid spacing $\frac{\Delta x}{2}$.
\end{remark}

If we only focused on the central Active Flux method, we might stop here since diagonal norm/mass matrices offer the best structure for many types of equations, e.g., in the presence of variable coefficients or curved grids \cite{svard2004coordinate}.
However, we are interested in the upwind version of the Active Flux method, as this version is more commonly used in practice. It will become clear in the following sections that we cannot use the diagonal mass matrix $M$ \eqref{eq:central-AF-diagonal-mass} for the upwind case. Therefore, here, we continue to explore a more general class of mass matrices.

Among diagonal matrices, \eqref{eq:central-AF-diagonal-mass} is unique up to a scalar factor\footnote{Since a direct proof of this fact would interrupt the flow of arguments here, we refer to Lemma~\ref{lem:central-AF-banded-mass}; requiring that $M$ be diagonal yields $m_{vv} = 0$ and $m_v = 3 m_p$, i.e., uniqueness of $M$ up to a scalar factor $m_p$.}.
Allowing for non-diagonal $M$, a mass matrix which induces the SBP property is not unique, as we can see in the following discussion. The additional parameters introduced here will be required for the upwind Active Flux method. Firstly, the following lemma shows the non-uniqueness of a symmetric mass matrix satisfying $ M D = -D^T M$.
\begin{lemma}
\label{lem:central-AF-banded-mass}
  All the symmetric, pentadiagonal matrices (up to boundary terms) which ensure that the derivative operator $D$ of the central
  semi-discrete Active Flux method~\eqref{eq:central-AF} is skew-symmetric
  with respect to $M$, are of the form
  \begin{equation}
  \label{eq:central-AF-banded-mass}
    M = \Delta x \left(\begin{array}{cc|cc|cc|cc|c}
      &  & \text{\rotatebox{10}{$\ddots$}} &  && \\
      & m_{vv} & \frac{m_v-3 m_p}{2} & m_v & \frac{m_v-3 m_p}{2} & m_{vv} \\\cline{1-8}
      && \frac{3 m_p - m_v+2m_{vv}}{6}& \frac{m_v-3 m_p}{2} & m_p & \frac{m_v-3 m_p}{2} & \frac{3 m_p - m_v+2m_{vv}}{6} & \\
      &&& m_{vv} & \frac{m_v-3 m_p}{2} & m_v & \frac{m_v-3 m_p}{2} & m_{vv} \\\cline{3-9}
      \multicolumn{4}{l|}{}  &&& \text{\rotatebox{10}{$\ddots$}} &  &  \\
    \end{array}\right)
  \end{equation}
  with $m_p, m_v, m_{vv} \in \R$.
\end{lemma}

\begin{proof}
Due to the alternate setting of point values and cell averages in the Active Flux method, we use the ansatz
  \begin{equation}
    M = \Delta x \left(\begin{array}{cc|cc|cc|cc|c}
       &  & \ddots & &&   \\
      & m_{vv} & m_{vp} & m_v & m_{vp}' & m_{vv} \\\cline{1-8}
      && m_{pp} & m_{vp}' & m_p & m_{vp} & m_{pp} &\\
      &&& m_{vv} &  m_{vp} & m_v & m_{vp}' & m_{vv} \\\cline{3-9}
      \multicolumn{4}{l|}{}&& & \ddots &  &  \\
    \end{array}\right),
  \end{equation}
  where $m_v$ and $m_p$ are on the diagonal of $M$ corresponding to the cell average (volume) and the point values, respectively. The skew symmetry  $M D + D^T M = 0$ requires
  \begin{align}
    m_{vp}' - m_{vp} &= 0,\\  3 m_{pp} + m_{vp} - m_{vv} &= 0, \label{eq:condM2} \\
    3 m_p - 3 m_{pp} + m_{vp}' - m_v + m_{vv} &= 0, \label{eq:condM3}
  \end{align}
  from which we directly obtain $m_{vp}' = m_{vp}$. Using this and adding the last two equations \eqref{eq:condM2}, \eqref{eq:condM3} we can solve for $m_{vp}$ and then for $m_{pp}$. Thus, we obtain
  \begin{align}\label{eq:mpp_mvp}
   m_{pp} &= \frac{3 m_p - m_v+2m_{vv}}{6}, \quad m_{vp} = \frac{m_v - 3 m_p}{2}\,,
  \end{align}
  which leaves us with the structure of $M$ as presented in \eqref{eq:central-AF-banded-mass}.
\end{proof}

\begin{remark}
    The banded mass matrix \eqref{eq:central-AF-banded-mass} reduces to the diagonal mass matrix in \eqref{eq:central-AF-diagonal-mass} for $m_{vv} = 0$, $m_p = m_v / 3$, and $m_v = \frac34$.
\end{remark}

Secondly, for the matrix $M$ in \eqref{eq:central-AF-banded-mass}
to qualify as a mass matrix of an SBP operator, we require $M$ to be
positive definite. In the following, we restrict the investigation
to $m_{vv} = 0$, since this will be required later for
the upwind version of Active Flux method. Moreover, setting
$m_{vv} \ne 0$ would result in a bigger stencil of $M$
in comparison to the stencil of the derivative operators. The following lemma determines a range of parameters for which $M$ is positive definite and gives one case of a positive semidefinite $M$ with multiplicity one of the corresponding zero eigenvalue.

\begin{lemma}
\label{lem:central-AF-banded-mass-definite}
  For $m_{vv} = 0$, the matrix $M$ in \eqref{eq:central-AF-banded-mass} has the following
  properties.
  \begin{enumerate}
    \item 
    If $2 m_v / 9 < m_p < 2 m_v / 3$ and $m_v > 0$, then $M$ is positive definite.
    \item If $m_p = 2 m_v / 9$
          and $m_v > 0$ then $M$ is positive semidefinite with an eigenvalue $0$ of multiplicity $1$.
  \end{enumerate}
\end{lemma}
The technical proof of Lemma~\ref{lem:central-AF-banded-mass-definite} involving properties of circulant matrices can be found in Appendix~\ref{sec:circulant-proof-of-lem:central-AF-banded-mass}.

For $m_{vv} = 0$ and $m_p \ne 2 m_v / 3$, it is possible to rescale the matrix $M$
in \eqref{eq:central-AF-banded-mass} to satisfy
the classical constraint $\vec{1}^T M \vec{1} = \xmax - \xmin$. Thus, we obtain
\begin{theorem}
\label{thm:central-SBP-operators}
    The central version \eqref{eq:findiffgeneralcentral}
    of the semi-discrete Active Flux method can be
    formulated using SBP operators with derivative matrix
    \begin{equation}
        D = \frac{1}{\Delta x} \left(\begin{array}{cc|cc|cc|c}
          \ddots &&&\\
          -1 & 0 & 1 &&\\\cline{1-6}
      1 & -3 & 0 & 3 & -1 &&\\
      && -1 & 0 & 1 &&\\\hline
      && 1 & -3 & 0 & 3 & -1 \\
      &&&& -1 & 0 & 1 \\\cline{3-7}
      \multicolumn{4}{l|}{}&&& \ddots
      \end{array}\right),
      \label{eq:D-central-AF}
    \end{equation}
    and positive definite mass matrix
    \begin{equation}
          M = \frac{3 \Delta x}{8 m_v - 12 m_p} \left(\begin{array}{cc|cc|cc|cc|c}
          &  & \text{\rotatebox{10}{$\ddots$}} &  && \\
          &  & \frac{m_v-3 m_p}{2} & m_v & \frac{m_v-3 m_p}{2} &  \\\cline{1-8}
          && \frac{3 m_p - m_v}{6}& \frac{m_v-3 m_p}{2} & m_p & \frac{m_v-3 m_p}{2} & \frac{3 m_p - m_v}{6} & \\
          &&& & \frac{m_v-3 m_p}{2} & m_v & \frac{m_v-3 m_p}{2} & \\\cline{3-9}
          \multicolumn{4}{l|}{} && & \text{\rotatebox{10}{$\ddots$}} &  &  \\
        \end{array}\right),
    \end{equation}
 for $m_v>0$ and $2 m_v / 9 < m_p < 2 m_v / 3$,  satisfying $\vec{1}^T M \vec{1} = \xmax - \xmin$.
\end{theorem}

Next, we will use the additional degrees of freedom introduced by the banded structure of the mass matrix $M$ in Theorem~\ref{thm:central-SBP-operators} to generalize the results to the upwind version of the Active Flux method.

\subsection{Upwind version of the Active Flux method}
\label{subsec:analysis_upwind}

We first collect the formulas for the upwind semi-discrete Active Flux method. The cell average update is
\begin{equation}
\label{eq:upwind-AF-average}
  \frac{\dif}{\dif t} u_i + \frac{u_{i+1/2} - u_{i-1/2}}{\Delta x} = 0,
\end{equation}
the upwind point value update
for positive advection velocity (see \eqref{eq:upwind-AF-Dm}) is
\begin{equation}
  \frac{\dif}{\dif t} u_{i+1/2} + \frac{2 u_{i-1/2} - 6 u_{i} + 4 u_{i+1/2}}{\Delta x} = 0,
\end{equation}
and for negative advection velocity (see \eqref{eq:upwind-AF-Dp}), we have
\begin{equation}
   \frac{\dif}{\dif t} u_{i+1/2} - \frac{-4 u_{i+1/2} + 6 u_{i+1} - 2 u_{i+3/2}}{\Delta x} = 0.
\end{equation}

The upwind versions of the Active Flux method correspond to the derivative
operators
\begin{equation}
\label{eq:upwind-AF}
  D_+ =
  \frac{1}{\Delta x} \left(\begin{array}{cc|cc|cc|c}
    \ddots &&&\\
    -1 & 0 & 1 &&\\\cline{1-6}
    && -4 & 6 & -2 && \\
    && -1 & 0 & 1 &&\\\cline{3-7}
    \multicolumn{4}{l|}{} & -4 & 6 & -2 \\
    \multicolumn{4}{l|}{} & -1 & 0 & 1 \\\cline{5-7}
    \multicolumn{6}{l|}{} & \ddots
  \end{array}\right),
  \quad
  D_- =
  \frac{1}{\Delta x}\left(\begin{array}{cc|cc|cc|c}
    \ddots &&& \\
    -1 & 0 & 1 &&\\\cline{1-6}
    2 & -6 & 4 &&&&\\
    && -1 & 0 & 1 &&\\\hline
    && 2 & -6 & 4 &&\\
    &&&& -1 & 0 & 1 \\\cline{3-7}
    \multicolumn{6}{l|}{}& \ddots
  \end{array}\right).
\end{equation}
Unfortunately, these operators are not SBP operators with the diagonal mass matrix
\eqref{eq:central-AF-diagonal-mass}. However we can still achieve this property using a banded
mass matrix as for the central approximation above. In addition, restricting the band width of the mass matrix to the band width of the derivative operators results in uniqueness of the admissible mass matrix $M$ up to a scalar factor.

\begin{lemma}
\label{lem:upwind-AF-banded-mass}
  The symmetric matrix
  \begin{equation}
  \label{eq:upwind-AF-banded-mass}
    M = \Delta x \left(\begin{array}{cc|cc|cc|cc|c}
       &  & \text{\rotatebox{10}{$\ddots$}} &  & & \\
      & 0 & -m_v / 2 & m_v & -m_v / 2 & 0 \\\cline{1-8}
      && m_v / 6 & -m_v / 2 & 2 m_v / 3 & -m_v / 2 & m_v / 6 & \\
      &&& 0 & -m_v / 2 & m_v & -m_v / 2 & 0 \\\cline{3-9}
      \multicolumn{4}{l|}{} && & \text{\rotatebox{10}{$\ddots$}} &  &  \\
    \end{array}\right)
  \end{equation}
  is --- up to a scalar factor $m_v \in \R$ --- the only symmetric matrix which is pentadiagonal (up to boundary terms) and ensures that the derivative operators $D_\pm$ of the upwind
  semi-discrete Active Flux method \eqref{eq:upwind-AF} are adjoint to each other
  with respect to $M$.
\end{lemma}
\begin{proof}
  We use the same ansatz as for the proof of Lemma \ref{lem:central-AF-banded-mass}.
  Then, the mutual adjointness
  $M D_+ + D_-^T M = 0$ requires $m_{vv}=0$ and
  \begin{align}
   m_{vp} - m_{vp}' &= 0
 &6 m_{pp} + 4 m_{vp}' + m_v &= 0 \\
 6 m_p + 4 m_{vp} - m_v +2 m_{vp}'&= 0
 &6 m_{pp} + 2 m_{vp} &= 0
  \end{align}
  This yields $m_{vp}=m_{vp}'=-m_v/2,\ m_p=2m_v/3,\ m_{pp}=m_v/6$ and hence \eqref{eq:upwind-AF-banded-mass}.
\end{proof}

\begin{remark}
    Note that the matrix $M$ in \eqref{eq:upwind-AF-banded-mass} is a special case of
the matrix in \eqref{eq:central-AF-banded-mass} with $m_p = 2 m_v / 3$.  In this case, $M$ is positive semidefinite as shown in the following lemma.
\end{remark}

\begin{lemma}
\label{lem:upwind-AF-banded-mass-pos-semidef}
    The matrix $M$ in \eqref{eq:upwind-AF-banded-mass} is positive semidefinite
    for $m_v > 0$ with an eigenvalue $0$ of multiplicity~$1$
    and eigenvector $\vec{1}$.
\end{lemma}
\begin{proof}
A simple calculation shows that the rows of $M$ sum up to zero. Thus, $\vec{1}$ is an eigenvector of $M$ with eigenvalue $0$. It is sufficient to consider the case $m_v=1$ to study its multiplicity and to determine the sign of the other eigenvalues. Using Lemma \ref{lem:circulant}, the matrix $B_k$ associated with $M$ for $m_v=1$ is given by (with $\theta = \frac{2\pi k}{n}$, $k=0, \dots, n-1$)
\begin{align}
B_k=\begin{pmatrix}
1 & -(1+e^{-i\theta})/2\\
-(1+e^{i\theta})/2 & (2+\cos{\theta})/3
\end{pmatrix}
\end{align}
and its characteristic polynomial is
\begin{align}
z^2 -  \frac{z}{3} (5 + \cos \theta) + \frac{1}{6} (1 -\cos \theta).
\end{align}
Its value at $z=0$ can only vanish if $\theta = 0$. Thus, this eigenvalue only appears once.
Furthermore, for $\theta\in(0,2\pi)$, the eigenvalues are positive since the roots of the characteristic polynomial are given by $z_{1,2}=\frac16\left(\cos\theta+5\pm\sqrt{(\cos\theta+5)^2 -6(1-\cos\theta)}\right)>0$.
\end{proof}

\begin{remark}
Classical finite difference SBP operators have the scaling $D \propto \Delta x^{-1}$ and $M \propto \Delta x$. Remark~\ref{rem:trapezoidalinterpretation} demonstrates that this expected scaling of the mass matrix $M$ also holds for the central Active Flux method. However, we cannot obtain the scaling $M \propto \Delta x$ by imposing $\vec{1}^T M \vec{1} = \xmax - \xmin$ since $M \vec{1} = \vec{0}$ for the upwind Active Flux method. Thus, we need to consider another way to fix the free parameter $m_v$ in \eqref{eq:upwind-AF-banded-mass}. To do so, we introduce the vector $\vec{1}_\text{avg} = (\dots, 0, 1, 0, 1, \dots)^T$ which has a value unity for all averages and zero for all point values. Then, we compute
\begin{align*}
	\frac{\vec{1}_\text{avg}^T M \vec{u}}{\Delta x \, m_v}
    = & \dots +\,\, u_{i-2} \,\,-\,\, u_{i-\frac32} \,\,+\,\, u_{i-1} \,\,-\,\, u_{i-\frac12} \,\,+\,\, u_i \,\,-\,\,u_{i+\frac12} \,\,+\,\, u_{i+1} \,\, - \dots
\end{align*}
In particular, we can choose $\vec{u} = \vec{1}_\text{avg}$, resulting in
\begin{align*}
	\frac{\vec{1}_\text{avg}^T M \vec{1}_\text{avg}}{\Delta x \, m_v}
    = & \dots +\,\, u_{i-2} \,\,+\,\, u_{i-1} \,\,+\,\, u_i \,\,+\,\, u_{i+1} \,\, + \dots
    =
    \frac{\xmax - \xmin}{\Delta x}.
\end{align*}
This shows that choosing $m_v = 1$ yields
$\vec{1}_\text{avg}^T M \vec{1}_\text{avg} = \xmax - \xmin$, i.e., the expected scaling of the mass matrix $M \propto \Delta x$ (when focusing only on the volume contributions).
\end{remark}

Relaxing the definiteness condition to allow for a positive semidefinite mass matrix $M$ still enables a stability proof for the Active Flux method via SBP properties as we will show in \cref{sec:stability}. At this point however, we still need to prove the upwind SBP property requiring that $M (D_+ - D_-)$ is negative semidefinite. This is the subject of the following lemma.

\begin{lemma}
\label{conj:upwind-AF-disipation}
    For $M$ in \eqref{eq:upwind-AF-banded-mass} and the upwind operators
    $D_\pm$ in \eqref{eq:upwind-AF}, the matrix $M (D_+ - D_-)$ is
    symmetric and negative semidefinite for $m_v > 0$.
\end{lemma}
\begin{proof}
  Concerning the symmetry, we have
  \begin{equation}
    \setcounter{MaxMatrixCols}{20}
    \frac{M (D_+ - D_-)}{m_v}
    =
    \left(\begin{array}{c|cc|cc|cc|cc|cc|cc|c}
         &&\text{\rotatebox{10}{$\ddots$}}&  &  &  &  &
         \\\cline{1-9}
         \dots & -8/3 & 5 & -6 & 5 & -8/3 & 1 & -1/3 &\\
         \dots & 1 & -3 & 5 & -6 & 5 & -3 & 1 & 0 &\\\cline{1-11}
        & -1/3 & 1 & -8/3 & 5 & -6 & 5 & -8/3 & 1 & -1/3 &\\
       &&0 & 1 & -3 & 5 & -6 & 5 & -3 & 1 & 0 &\\\cline{2-13}
     \multicolumn{3}{l|}{} & -1/3 & 1 & -8/3 & 5 & -6 & 5 & -8/3 & 1 & -1/3 &\\
      \multicolumn{3}{l|}{}&& 0 & 1 & -3 & 5 & -6 & 5 & -3 & 1 & 0& \\\cline{4-13}
      \multicolumn{5}{l|}{}&&&&& \text{\rotatebox{7}{$\ddots$}}&  &  &  &   \\
    \end{array}\right).
  \end{equation}

 Finally, using Lemma \ref{lem:circulant}, one finds, calling $\theta := \frac{2\pi k}{n}$, $k=0,\dots,n-1$, the eigenvalues $\lambda(\theta)$ to be 0 ($n$ times) and 
 \begin{align}
   \lambda(\theta) = -\frac23 \left(18 + 17 \cos \theta + \cos\left(2\theta \right)\right) < 0, \qquad \forall \theta.
 \end{align}
\end{proof}

\begin{remark}
  Following the nomenclature introduced in \cite{hicken2024constructing}
  for central SBP operators, we combine the previous results
  (\Cref{lem:upwind-AF-banded-mass},
  \Cref{lem:upwind-AF-banded-mass-pos-semidef}, and \Cref{conj:upwind-AF-disipation})
  and say that the semi-discrete Active Flux upwind operators $D_\pm$
  are \emph{degenerate} upwind SBP operators when combined with the
  semidefinite mass matrix $M$ \eqref{eq:upwind-AF-banded-mass}.
\end{remark}

\section{Stability of the Active Flux method via SBP properties}\label{sec:stability}

Since the central version of the Active Flux method
can be formulated using
SBP operators, we directly obtain stability.
\begin{corollary}
    The central version \eqref{eq:findiffgeneralcentral}
    of the semi-discrete Active Flux method is stable
    for the linear advection \eqref{eq:linad} with periodic boundary conditions.
\end{corollary}
\begin{proof}
    This follows immediately from
    Theorem~\ref{thm:central-SBP-operators}
    and classical stability properties of SBP operators.
\end{proof}

Unfortunately, the matrix $M$ for which the upwind Active Flux operators $D_\pm$
are mutually adjoint is only positive semidefinite. Thus, the classical energy
stability proof does not guarantee stability immediately. However, based on
\Cref{lem:upwind-AF-banded-mass-pos-semidef}, we can still obtain a stability result
for the Active Flux method using the framework of SBP operators.

The idea is as follows. We show below that the states associated to zero energy are uniform (in space) constants. Thus, if the energy is decaying, it is, generally speaking, possible that a uniform constant could grow: our mass matrix has constants in its kernel and its energy is zero.  However, the time evolution of uniform constants is easy to check separately, and we find that uniform constants do not grow in time.

\begin{theorem}
  The upwind Active Flux semi-discretization
  \begin{equation}
     \frac{\dif}{\dif t} \vec{u} + D_- \vec{u} = \vec{0}
  \end{equation}
  using \eqref{eq:upwind-AF-average} and
  \eqref{eq:upwind-AF-Dm} of the linear advection equation \eqref{eq:linad} with $a = 1$
  with periodic boundary conditions is stable.
\end{theorem}
\begin{proof}
  Using \Cref{lem:upwind-AF-banded-mass-pos-semidef}, we decompose the vector
  $\vec{u} = \vec{u}_0 + \vec{u}_\perp$ into the components
  $\vec{u}_0 \propto \vec{1}$ in the kernel of $M$ and
  $\vec{u}_\perp \perp \vec{1}$ in its orthogonal complement.
  Following the usual upwind SBP argumentation
  \eqref{eq:upwind-SBP-stability}, we obtain
  \begin{equation}
  \begin{aligned}
    \frac{\dif}{\dif t} \| \vec{u} \|_M^2
    &=
    2 \vec{u}^T M  \frac{\dif}{\dif t} \vec{u}
    =
    -2 \vec{u}^T M D_- \vec{u}
    =
    - \vec{u}^T M D_- \vec{u} + \vec{u}^T D_+^T M \vec{u}
    \le
    0,
  \end{aligned}
  \end{equation}
  where we have used \Cref{conj:upwind-AF-disipation}.
  Since $M \vec{1} = \vec{0}$,
  $\| \vec{u} \|_M = \| \vec{u}_\perp \|_M$ and the component in
  the orthogonal complement is stable. The component in the kernel
  of $M$ is constant in time and thus also stable, since
  \begin{equation}
     \frac{\dif}{\dif t} \vec{u}_0 = -D_- \vec{u}_0 = \vec{0}.
  \end{equation}
  Thus, the method is stable.
\end{proof}
The Theorem and Proof above hold analogously for the linear advection equation \eqref{eq:upwind-AF-Dp} with the opposite speed when using $D_+$.

\section{Numerical experiment}\label{sec:experiments}

To exemplify the theoretical stability results, we perform a numerical
experiment based on the linear advection equation
$\partial_t u + \partial_x u = 0$
with periodic boundary conditions in the interval $[0, 2\pi]$ and
the initial condition $u(0, x) = \exp(\sin(x))$. We use the
central and upwind semi-discrete Active Flux methods with
10 and 50 volumes, respectively,
and integrate the resulting ODE system in time using the five-stage,
third-order explicit Runge-Kutta method
RK3(2)5\textsubscript{F}[3S\textsuperscript{*}\textsubscript{+}]
of \cite{ranocha2021optimized}
with time step size $\Delta t = \Delta x / 2$. We chose this method
since it has been optimized for hyperbolic problems and includes an
interval on the imaginary axis in its stability interval as required
for stability of central operators.

The energy analysis presented in the previous sections concentrates on the effect of the spatial semi-discretization.
When integrating the resulting ODE in time, the energy behavior also depends on the time integration scheme.
While many explicit Runge-Kutta methods can be shown to add additional dissipation for linear problems \cite{tadmor2002semidiscrete,sun2017stability,ranocha2018L2stability,sun2019strong}. As a remark, we mention here that the analysis of nonlinear problems is more complex \cite{lozano2018entropy,ranocha2021strong,ranocha2020energy}.
Since we are interested in a quadratic energy, we could use a symplectic Runge-Kutta method such as the implicit midpoint rule to preserve the energy evolution caused by the semi-discretization in space, cf.\ sections~IV.2.1 and IV.4.1 of \cite{hairer2006geometric}.
However, we would like to avoid introducing fully implicit methods for hyperbolic problems.
Thus, we apply the relaxation approach discussed and analyzed in \cite{ketcheson2019relaxation,ranocha2020relaxation,ranocha2020general} to modify an explicit Runge-Kutta method slightly so that the fully discrete energy evolution is purely caused by the spatial semi-discretization without additional influence of the time integration method.
To demonstrate the dissipative effect of the time integration method, we also show results without relaxation in time.

\begin{figure}[htb]
  \centering
  \includegraphics[width=\textwidth]{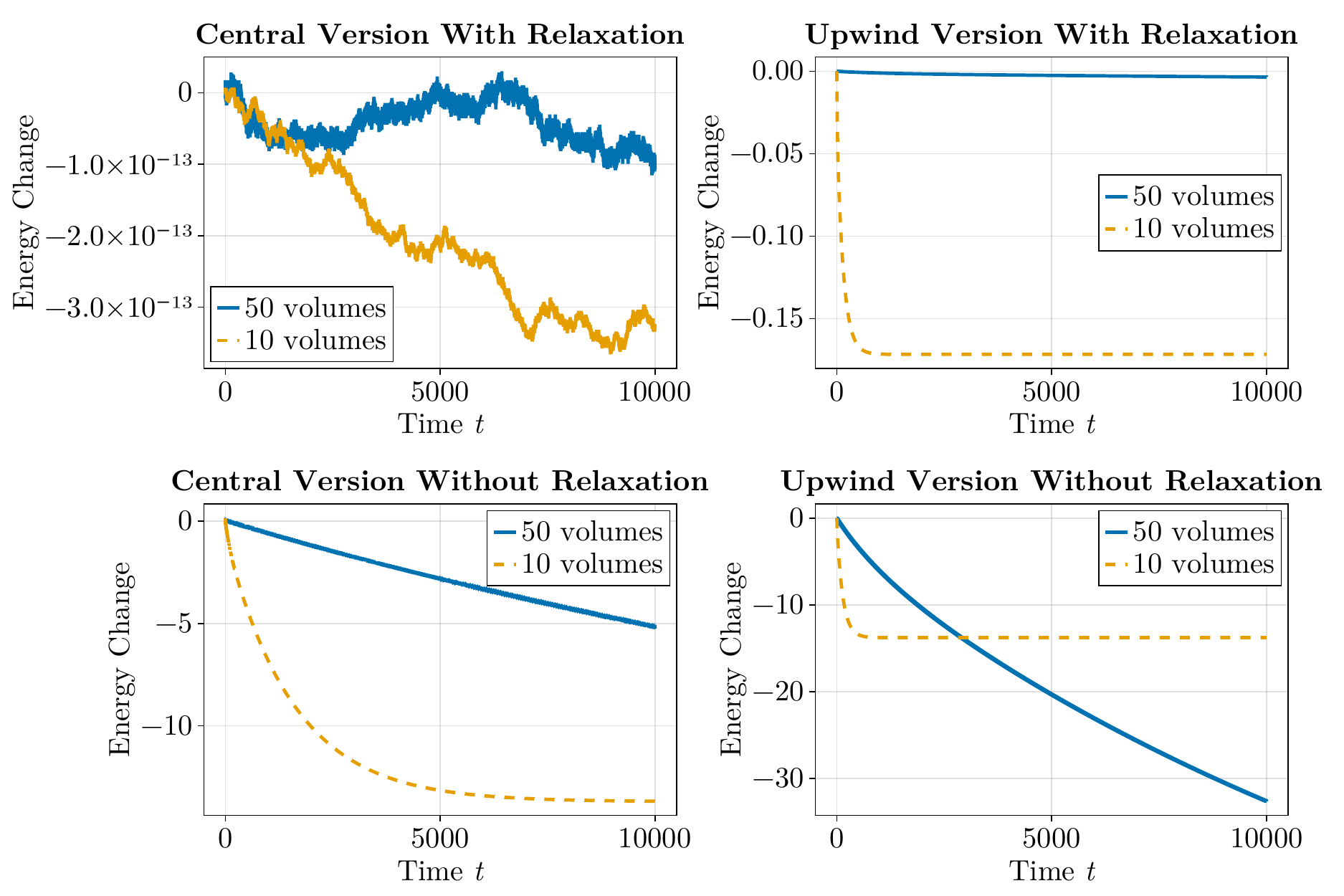}
  \caption{Change of the discrete energy
           $\| \vec{u} \|_M^2 = \pmb{u}^T M \pmb{u}$
           compared to the initial energy. For the central version,
           we use the diagonal mass matrix
           \eqref{eq:central-AF-diagonal-mass}; for the upwind version,
           we use the pentadiagonal (up to boundary terms) mass matrix
           \eqref{eq:upwind-AF-banded-mass} with $m_v = 1$.}
  \label{fig:energies}
\end{figure}

The results are shown in Figure~\ref{fig:energies}. In accordance with the
theory, the energy is conserved up to roundoff errors for the central
version while the upwind version results in energy dissipation if relaxation in time is used.
Without relaxation, the explicit time integration method introduces additional dissipation.
All Julia \cite{bezanson2017julia} code and instructions required to reproduce these numerical experiments are available in our reproducibility repository \cite{barsukow2025stabilityRepro}.
We used CairoMakie.jl \cite{danisch2021makie} to create the figures.

\section{Nullspace consistency of the Active Flux upwind SBP operator}
\label{sec:nullspace consistency}

An additional property of the Active Flux difference operators which is not directly related to their upwind SBP property is nullspace consistency. This property is interesting on its own merit since it transfers an additional algebraic property from the continuous to the discrete setting. Nullspace consistency is a rather recently studied mimetic property in the context of derivative operators and numerical methods for PDEs. It signifies that the nullspace of the continuous operator is correctly transferred to its discrete counterpart. In the context of SBP schemes, this concept was first introduced in \cite{svard2019convergence}, connected to the investigation of convergence rates of FD schemes. According to the definition of nullspace consistency, for a nullspace consistent finite difference operator $D_k$ approximating the continuous operator $\partial_x^k$, the kernel of $D_k$ can be mapped one-to-one to the kernel of $\partial_x^k$. 

In the context of periodic problems such as those studied in this work, nullspace consistency prevents spurious modes, that are contained in the discrete kernel of $D_k$ but not in its continuous counterpart, from persisting in the numerical solution for all times. In the absence of round-off and cancellation errors, such spurious modes are unaffected by the numerical scheme and remain constant in time, since the discrete derivative operator maps them to zero. A spurious mode may hence exist in the initial data, e.g. in an underresolved simulation, or be created due to round-off errors. The continuous linear advection equation transports this mode at the speed determined by the advection coefficient, while the semi-discrete scheme leaves it unchanged in time. Since we have a linear problem, due to energy stability, such modes will not get amplified since their energy cannot grow. However, if the scheme is not nullspace consistent, there is also no dissipation mechanism which can damp out this spurious high-frequency mode in the periodic setting. In contrast, nullspace consistency allows for such a damping mechanism. The following lemma shows that upwind Active Flux difference operators are nullspace consistent.

\begin{lemma}
  The Active Flux difference operators $D_-$ and $D_+$ are nullspace consistent, i.e.,
  $D_\pm \vec{u} = \vec{0} \iff \vec{u} \propto \vec{1}$.
\end{lemma}

\begin{proof}
Due to symmetry, it is sufficient to consider $D_-$ as the proof for $D_+$ follows analogously.
The nullspace of the linear continuous operator $\partial_x$ consists of the functions constant in $x$. Thus, we require that the only vectors $\vec{u}$ satisfying $D_- \vec{u} = \vec{0}$ are the constants $\vec{u} = c \vec{1}$.

Obviously, all vectors of the form $\vec{u} = c \vec{1}$ satisfy $D_- \vec{u} = \vec{0}$. It remains to show that all the vectors $\vec{u}$ with $D_- \vec{u} = \vec{0}$ are of the form $\vec{u} = c \vec{1}$. Then, all vectors in the nullspace of $D_-$ have a counterpart in the nullspace of $\partial_x$.
The linear difference scheme resulting from the condition $D_- \vec{u} = \vec{0}$ is
\begin{align*}
 \forall i&\colon& -u_{i-\frac12}+u_{i+\frac12} &=0,\\
 \forall i&\colon& 2u_{i-\frac12} - 6 u_i+4u_{i+\frac12} &=0.
\end{align*}
From the first equation, we obtain the existence of a constant $c\in\R$ such that $u_{i+\frac12} = c$ is constant for all $i$. Inserting this into the second equation then yields $u_i = c$ for all $i$, i.e., $\pmb{u} = c \pmb{1}$.
\end{proof}

The central version of the Active Flux method is not nullspace consistent. This is in line with the lack of nullspace consistency for periodic problems of the classical central second-order finite difference operator approximating the first derivative based on the approximation $\partial_x u(x_i) \approx \frac{u_{i+1}-u_{i-1}}{2\Delta x}$.

\begin{lemma}
  The central Active Flux difference operator $D$ in \eqref{eq:central-AF} has a two-dimensional nullspace spanned by the vectors $\vec{1}$ and $(1,-1,1,-1,\ldots,1,-1)^T$. Therefore, the central version \eqref{eq:findiffgeneralcentral} of the Active Flux method is not nullspace consistent.
\end{lemma}

\begin{proof}
As before, all vectors of the form $\vec{u} = c \vec{1}$ obviously satisfy $D \vec{u} = \vec{0}$. Furthermore, the linear difference scheme resulting from the condition $D \vec{u} = \vec{0}$ is given by
\begin{align*}
 \forall i&\colon& -u_{i-\frac12}+u_{i+\frac12} &=0,\\
 \forall i&\colon& u_{i-\frac12} - 3 u_i + 3u_{i+1} - u_{i+\frac32} &=0.
\end{align*}
From the first equation, we again obtain that $u_{i+\frac12} = c$ is constant for all $i$. Inserting this into the second equation now yields $u_{i+1} = u_{i}=d$ for all $i$, for a potentially different constant $d$.
Setting $c=1$ and $d=\pm1$ we obtain the two linearly independent vectors $\vec{1}$ and $(1,-1,1,-1,\ldots,1,-1)^T$ which span the nullspace of $D$.
\end{proof}

\begin{remark}
On finite non-periodic domains, nullspace consistency for central difference approximations can be restored by suitable boundary closures \cite{svard2019convergence}.
\end{remark}

\begin{remark}
The lack of nullspace consistency is not specific to the central Active Flux method, but just a consequence of its central character.
Indeed, it is well-known that any periodic SBP operator with an even number of degrees of freedom is not nullspace consistent:
A periodic SBP operator is skew symmetric (with respect to the inner product induced by the mass matrix).
Thus, it is unitarily diagonalizable, all eigenvalues are on the imaginary axis, and eigenvalues are either real or come in complex conjugate pairs.
Zero is an eigenvalue because of consistency of the SBP operator (mapping constants to zero).
Thus, if a periodic SBP operator has an even number of degrees of freedom, it needs to have an eigenspace associated to the eigenvalue zero of dimension at least two.
For example, all central finite difference operators with an even number of nodes have grid oscillations $(1,-1,1,-1,\ldots,1,-1)^T$ in their nullspace.
Fourier collocation methods map such grid oscillations (the Nyquist frequency) to zero as well.
\end{remark}

\section{Summary and discussion}

To the best of our knowledge, we have proven stability using the summation-by-parts technique for an Active Flux method for the first time. To this end, we considered the linear advection equation on an interval with periodic boundary conditions. The Active Flux method, a finite volume method based on a globally continuous approach, combines cell average updates with point updates at the cell edges. We have studied two ways of updating the point values.
First the point value updates using a \emph{central discretization} \eqref{eq:central_AF_point_values} are chosen for purposes of theoretical analysis in \cref{subsec:analysis_central}. In \cref{subsec:analysis_upwind} the commonly used \emph{upwind} (see \eqref{eq:upwind-AF-Dm} and \eqref{eq:upwind-AF-Dp}) version of the point value update is considered. Stability is shown for both cases in \cref{sec:stability}.  In Section~\ref{sec:experiments}, we provide numerical examples illustrating stability, consistent with the theoretical results in this paper. In Section~ \ref{sec:nullspace consistency}, we showed that the upwind Active Flux method is nullspace consistent while the central one is not. As an additional bonus, the relation between energy stability via summation by parts and von Neumann stability is shown in \cref{sec:von-Neumann-stability} for our PDE at hand.

While the central Active Flux version can be formulated using SBP operators with diagonal and positive definite mass matrix, the upwind version requires non-diagonal 
positive semi-definite mass matrices.
This does not impact the stability of the method for the test problem under consideration. However, non-diagonal mass matrices are troublesome for variable coefficients \cite{Nordström2006375}, curved grids \cite{svard2004coordinate}, and nonlinear problems \cite{Fisher2013353}, since they do not commute with non-constant diagonal matrices.
In special cases such issues 
can sometimes be solved by adding sufficient stabilization/dissipation 
\cite{mattsson2013solution}.
It is currently unclear whether the dissipation introduced by the upwind Active Flux operators is sufficient to stabilize such situations.
Furthermore, we 
stress that our findings have to be generalized to the multi-dimensional case. 
Indeed, the Active Flux method is not used in a tensor product fashion, but is genuinely multi-dimensional, i.e., it uses 
different operators that do not follow directly from the 1D operators. We currently do not know how 
mass matrices might be structured in the multi-dimensional case and
see this as an opportunity for interesting future research.
The apparent success of the Active Flux method suggests that there may be even more hidden structures and properties that we are currently unaware of.

\section*{Acknowledgments}

The authors acknowledge funding by the Deutsche Forschungsgemeinschaft (DFG, German Research Foundation) within \emph{SPP 2410 Hyperbolic Balance Laws in Fluid Mechanics: Complexity, Scales, Randomness (CoScaRa)}, project numbers 525941602 (W.B., C.K., L.L.),  526073189 (S.O.) and 526031774 (H.R.). 
H.R. additionally acknowledges the DFG individual research grant 513301895. J.N. was supported by Vetenskapsrådet, Sweden [no. 2021-05484 VR] and University of Johannesburg Global Excellence and
Stature Initiative Funding.

\appendix

\section{Properties of block circulant matrices and Proof of Lemma~\ref{lem:central-AF-banded-mass-definite}}
\label{sec:circulant-proof-of-lem:central-AF-banded-mass}

Block circulant matrices \cite{trapp1973blockcirculant,davis1979} are an extension of circulant matrices and possess the general form
\begin{equation}\label{eq:blockcirc}
A = \begin{pmatrix}
       A_0 & A_1 & & & A_{n-1}\\
       A_{n-1} & A_0 & A_1 & & \\
       & A_{n-1} & A_0 & A_1 & &\\
       & & & \ddots & &\\
      A_1 & & &A_{n-1} & A_0
      \end{pmatrix}\in \R^{nm\times nm},
\end{equation}
with $A_k\in\R^{m\times m}$ for $k=0,\ldots,n-1$.

The semi-discrete Active Flux method in one space dimension on a periodic domain yields block circulant matrices $D_+,D_-$ and $M$ and consequently also $M(D_+-D_-)$ has a block circulant structure. Due to the alternating difference formulas for point values and for cell averages, the involved matrices hereby possess a specific structure of $2\times2$ subblocks
\[A_k=\begin{pmatrix}
a_k & a_{k+1}\\
b_k & b_{k+1}\end{pmatrix},\quad k=0,\ldots,n-1.\]
We may thus exploit the theory of block circulant matrices to obtain eigenvalues and eigenvectors of the given matrices.

\begin{lemma}[Eigenvalues of block circulant matrices] \label{lem:circulant}
The eigenvalues of the $nm \times nm$ matrix \eqref{eq:blockcirc} are the eigenvalues of all matrices
\begin{equation}\label{eq:Bk}
B_k = A_0+r^k A_1+ \ldots r^{k(n-1)}A_{n-1}, \qquad k=0,1,\ldots,n-1,
\end{equation}
with $r=\mathrm{e}^{2\pi \ii/n}$ the roots of unity. If $\lambda$ is an eigenvalue of $B_k$ with corresponding eigenvector $v$, then $\lambda$ is also an eigenvalue of $A$ with eigenvector $(v^\mathrm{T}, r^k v^\mathrm{T}, r^{2k} v^\mathrm{T}, \ldots, r^{(n-1)j} v^\mathrm{T} )^\mathrm{T}$.
\end{lemma}
\begin{proof}
 Consider the following ansatz for an eigenvector of $A$
\begin{align}
 e_k = (v^\mathrm{T}, r^k v^\mathrm{T}, r^{2k} v^\mathrm{T}, \ldots, r^{(n-1)j} v^\mathrm{T} )^\mathrm{T}, \qquad k = 0, 1, \ldots, n-1,
\end{align}
for some yet to be determined $v \in \mathbb C^m$. Then,
$A e_k = \lambda_k e_k$ if
\begin{align}
 A e_k &= \begin{pmatrix}
    A_0 v + r^k A_1 v + \ldots + r^{(n-2)k}A_{n-2} v + r^{(n-1)k} A_{n-1} v \\
    A_{n-1} v + r^k A_0 v + r^{2k} A_1 v + \ldots + r^{(n-1)k} A_{n-2} v\\
    \vdots \\
    A_1 v + \ldots + r^{(n-3)k} A_{n-2} v + r^{(n-2)k} A_{n-1} v + r^{(n-1)k} A_0 v
  \end{pmatrix} = \lambda_k \begin{pmatrix}
                v \\ r^k v \\ \vdots \\ r^{(n-1)k} v
               \end{pmatrix}.
\end{align}
The $i$-th component of this equation is ($i = 0, \ldots, n-1$)
\begin{align}
\sum_{j = n-i}^{n-1}  r^{(j-n+i)k}  A_j v  + r^{ik} \sum_{j=0}^{n-1 - i} r^{jk} A_j v  &= \lambda_k r^{ik} v, \label{eq:singlelinecirculant}
\end{align}
which reduces to just one equation for all $i$ upon the choice $r = \exp\left( \frac{2\pi \ii}{n} \right )$.
Then, the eigenvalue problem reduces to
\begin{align}
\sum_{j=0}^{n-1} r^{jk} A_j v = B_k v &= \lambda_k v.
\end{align}
\end{proof}

As shown in \cite{trapp1973blockcirculant}, all eigenvalues and eigenvectors of $A$ can be obtained in this manner.
In addition, while all circulant matrices can be diagonalized, block circulant matrices allow for block diagonalization. A block circulant matrix does not necessarily possess a complete set of linearly independent eigenvectors; however, it does so if the matrices $B_k$ can all be diagonalized. In this case, the block circulant matrix can be diagonalized as well. 

\begin{lemma}\label{lem:blockdiag}
Let $F\in\C^{n\times n}$ be the unitary matrix with entries $F_{jk}=\frac1{\sqrt{n}}r^{(j-1)(k-1)},\ r=e^{2\pi \ii/n}$, and $I_m=\operatorname{diag}(1,\ldots,1)\in \R^{m\times m}$.
For a block circulant matrix $A\in\R^{mn\times mn}$ as in \eqref{eq:blockcirc}, we have
\begin{equation}\label{eq:blockdiag}
(F^\ast\otimes I_m) A (F\otimes I_m) = \operatorname{diag}(B_0,\ldots,B_{n-1}),
\end{equation}
with $B_k,\ k=0,\ldots,n-1$ given in \eqref{eq:Bk}.
\end{lemma}
\begin{proof}
Using $$F\otimes I_m=\frac{1}{\sqrt{n}}\begin{pmatrix}I_m & I_m & I_m &\cdots & I_m\\
I_m & rI_m & r^2I_m & \cdots & r^{n-1}I_m\\
I_m & r^2I_m & r^4I_m & \cdots & r^{2(n-1)}I_m\\
\vdots & & \vdots\\
I_m & r^{n-1}I_m & r^{2(n-1)}I_m & \cdots & r^{(n-1)^2}I_m\end{pmatrix}$$ and $r^n=1$, we first obtain
\[ A (F\otimes I_m) = \frac{1}{\sqrt{n}}\begin{pmatrix}
B_0 & B_1 & B_2 & \cdots & B_{n-1}\\
B_0 & B_1r & B_2r^2 & \cdots & B_{n-1}r^{n-1}\\
\vdots & &&&\vdots\\B_0 & B_1r^{n-1} & B_2r^{2(n-1)} & \ldots & B_{n-1}r^{(n-1)^2}
\end{pmatrix}
\]
with $B_k$ given in \eqref{eq:Bk}.
Using $F^\ast_{jk}=\frac1{\sqrt{n}}\bar{r}^{(j-1)(k-1)}$ with $\bar{r}=e^{-2\pi \ii/n}$, i.e.,
\[F^\ast\otimes I_m=\frac{1}{\sqrt{n}}\begin{pmatrix}I_m & I_m & I_m &\cdots & I_m\\
I_m & \bar{r}I_m & \bar{r}^2I_m & \cdots & \bar{r}^{n-1}I_m\\
I_m & \bar{r}^2I_m & \bar{r}^4I_m & \cdots & \bar{r}^{2(n-1)}I_m\\
\vdots & & \vdots\\
I_m & \bar{r}^{n-1}I_m & \bar{r}^{2(n-1)}I_m & \cdots & \bar{r}^{(n-1)^2}I_m\end{pmatrix},\]
we have
\[(F^\ast\otimes I_m) A (F\otimes I_m)=\frac1n\begin{pmatrix}
nB_0 & p(r)B_1 & p(r^2)B_2 & \cdots & p(r^{n-1})B_{n-1}\\
p(r^{n-1})B_0 & nB_1 & p(r)B_2 & \cdots& p(r^{n-2})B_{n-1}\\
\vdots & & & \ldots & \vdots\\
p(r)B_0 & p(r^2)B_1 & p(r^3)B_2 &\cdots & nB_{n-1}
\end{pmatrix},
\]
with $p(z)=1+z+z^2+\ldots+z^{n-1}$, where we exploit $\bar{r}^k=r^{n-k}$. Finally, for $z=r^k$ with $k\not=0$ we have $p(z)=0$, since the values $r^k,\ k=1,\ldots,n-1$ are the roots of $z^n=1$ which differ from $z=1$ and we may factorize $z^n-1=(z-1)(1+z+z^2+\ldots +z^{n-1})$. Thus, the block-diagonal form \eqref{eq:blockdiag} is proven.
\end{proof}

\begin{proof}[Proof of Lemma~\ref{lem:central-AF-banded-mass-definite}]
We use Lemma~\ref{lem:circulant} and the notation described there. The eigenvalues of $M$ can then be determined by the eigenvalues of the matrices $B_k$ which in this case are given by
\begin{align}
 B_k = \left(\begin{array}{cc}m_p + 2m_{pp}\cos\theta & m_{vp}\left(1+e^{-\ii\theta}\right)\\m_{vp}\left(1+e^{\ii\theta}\right)&m_v\end{array}\right), \quad \theta := \frac{2\pi k}{n},
\end{align}
with $m_{pp}$ and $m_{vp}$ determined by \eqref{eq:mpp_mvp}. The characteristic polynomial $p$ of $B_k$ is given by
\begin{align}
 p(z)= z^2-tz+d
\end{align}
with $t$ and $d$ the trace and determinant of $B_k$, respectively,
\begin{align}
 t &= \mathrm{trace}(B_k) = \frac{1}{3} \left(3(m_p + m_v) + (3m_p-m_v) \cos\theta\right),\\
 d &= \det{B_k} =  \frac16 \left( -3(9 m_p^2 - 8 m_p m_v+m_v^2) - (9 m_p - 5 m_v)(3 m_p-m_v)\cos\theta   \right).
\end{align}
We are now interested in a range of parameters $m_v,m_p$ for which $B_k$ has only positive eigenvalues for any value of $\theta$.
From the solution formula of quadratic equations it is clear that for positive solutions one requires $t>0$ and $d>0$ which yields the following inequalities regarding the coefficients of the characteristic polynomial:
\begin{align}
3(m_p + m_v) + (3m_p-m_v) \cos\theta &>0, \label{eq:intermdineq1}\\
 -3(9 m_p^2 - 8 m_p m_v+m_v^2) - (9 m_p - 5 m_v)(3 m_p-m_v)\cos\theta   &>0. \label{eq:intermdineq2}
\end{align}
The latter inequality can be rewritten as
\begin{align}
 m_v^2 -(9 m_p - 2 m_v)(3m_p - 2m_v) - \left( (9 m_p - 2 m_v)(3 m_p-2m_v) + m_v^2\right ) \cos\theta   &>0 \label{eq:quadraticintermed}
\end{align}
i.e.
\begin{align}
 m_v^2 -s- \left( s + m_v^2\right ) \cos\theta   &>0
\end{align}
with $s:= (9 m_p - 2 m_v)(3 m_p-2m_v)$. For $m_v > 0$, we obtain $s<0$ precisely for $2 m_v / 9 < m_p < 2 m_v / 3$. In that case, by the strict triangle inequality, we have
\begin{align}
|s +  m_v^2| <  |s| + m_v^2 = m_v^2 - s.
\end{align}
Since the condition $a + b \cos \theta > 0$ for all $\theta$ is equivalent to $|b| < a$, one thus has proven \eqref{eq:intermdineq2} for the range of parameters considered in the first statement of the Lemma.
Furthermore, simply with $m_p, m_v > 0$ one obtains
\begin{align}
|3m_p - m_v| \leq 3|m_p| + |m_v| < 3 (m_p + m_v)
\end{align}
and thus \eqref{eq:intermdineq1} by the same argument. Thereby, we have proven the assertion that $M$ possesses only positive eigenvalues if $m_v>0$ and $2 m_v / 9 < m_p < 2 m_v / 3$.

Finally considering the case of positive semidefinite matrices $M$, we note that the  characteristic polynomial $p$ may have a zero eigenvalue, which is the case if and only if $d=\det{B_k}=0$. For instance, this happens if
\begin{itemize}
\item $\theta = 0$ and either $m_p = \frac{2 m_v}{9}$ or $m_p = \frac{2 m_v}{3}$,
\item $\theta = \pm \frac{\pi}{2}$ and $m_p = \frac{4 \pm \sqrt{7}}{9} m_v = \begin{cases} 0.738... m_v  > \frac23 m_v\,, \\ 0.150... m_v < \frac29 m_v\,. \end{cases}$
\end{itemize}
Analogously, other values of $\theta$ will likewise yield parameter values with $B_k$ possessing a zero eigenvalue but will not be further studied here.

Considering the case $m_p = \frac{2 m_v}{9}$, we have $\det B_k = \frac16 m_v^2(1 - \cos \theta)$, which only vanishes for $\theta=0$, and since we have $\mathrm{trace}(B_k) \not=0$ for $\theta=0,\ m_p = \frac{2 m_v}{9}$, the eigenvalue 0 has multiplicity 1. This proves the second statement of this lemma. The case $m_p = \frac{2 m_v}{3}$ is treated in Lemma \ref{lem:upwind-AF-banded-mass}.
\end{proof}

\section{Relation to von Neumann stability}
\label{sec:von-Neumann-stability}

Von Neumann stability uses the Fourier transform to make statements about the behavior in time of the $L^2$ norm of the solution. In its restricted setting of linear problems on periodic domains it can also be used to derive a mass matrix that satisfies the periodic SBP property, as will be shown now.

\subsection{General theory}

Since in the setting of the Active Flux method, cell averages and point values are independent degrees of freedom, and since their nature and their update equations differ, when performing von Neumann stability analysis they need to be associated to independent Fourier modes. We thus write
\begin{align}
u_j &= \sum_{\omega} \hat u^\text{avg}(\omega) \exp(\ii \omega j \Delta x) &
u_{j-\frac12} &= \sum_{\omega} \hat u^\text{point}(\omega) \exp(\ii \omega j \Delta x)
\end{align}
The summation is over all wave numbers $\omega$ that are compatible with the boundary conditions, i.e.,
\begin{align}
q_{j+n} = \sum_{\omega} \hat q^\text{avg}(\omega) \exp(\ii \omega (j+n) \Delta x)  \equiv \sum_{\omega} \hat q^\text{avg}(\omega) \exp(\ii \omega j \Delta x) = q_j
\end{align}
such that
\begin{align}
\omega \Delta x = 2 \pi \frac{k}{n} \qquad k = 0, \ldots, n-1.
\end{align}

Since the method is linear, Fourier modes do not mix. From now on, we thus perform all calculations for just one of them and $\omega$ is treated as a parameter. It is useful to introduce the translation operator
\begin{align}
\tau := \exp(\ii \omega \Delta x) = \exp\left(2 \pi \frac{k}{n}\ii \right).
\end{align}
Indeed, for a single Fourier mode,
\begin{align}
u_{j+1} &= \hat u^\text{avg}(\omega) \exp(\ii \omega (j+1) \Delta x) = \tau u_j &
\text{and similarly \quad } u_{j+\frac12} &= \tau u_{j-\frac12}.
\end{align}
The dependence of $\tau$ on $k$ shall not be made explicit in the notation.

In fact, $\tau = r^k$ that has appeared previously in Lemma \ref{lem:circulant}. The circulant matrices discussed in \cref{sec:circulant-proof-of-lem:central-AF-banded-mass} have $2\times 2$ blocks associated with one point value and one average. With $A$ defined in \eqref{eq:blockcirc}, inserting into $A \vec{u}$ a Fourier mode yields
\begin{align}
(A \vec{u})_i = \sum_{j = n-i}^{n-1} A_{j} \hat u \tau^{j-n+i} + \sum_{j = 0}^{n-i-1} A_j \hat u \tau^{i+j} \overset{\tau^n = 1}{=} \tau^i \sum_{j = 0}^{n-1} A_j \hat u \tau^{j} = \tau^i B_k \hat u
\end{align}
where $\hat u = \left( \begin{array}{c} \hat u^\text{point}\\\hat u^\text{avg} \end{array} \right)$. The parallel to Equation \eqref{eq:singlelinecirculant} is obvious. While the matrices $B_k$ have been so far only used to determine the eigenvalues of $A$, in fact it makes sense to call $B_k$ the \emph{Fourier symbol} $\hat A$ of $A$ that will be henceforth denoted by a hat.

Consider a generic numerical method
\begin{align}
 \frac{\dif}{\dif t} \vec u + A \vec u = 0.
\end{align}
Von Neumann stability analysis amounts to studying the eigenvalues of $\hat A$. Using the diagonalization $\hat R^{-1} \hat\Lambda \hat R$ of $\hat A$ in
\begin{align}
 \frac{\dif}{\dif t} \hat R \hat u + \hat \Lambda \hat R \hat u = 0
\end{align}
and calling $\hat U := \hat R \hat u$, each component $\hat U_\ell$ of $\hat U$ fulfills the ODE
\begin{align}
\frac{\dif}{\dif t} \hat U_\ell + \lambda_\ell \hat U_\ell = 0
\end{align}
solved by $\hat U_{\ell,0} \exp(-\lambda t)$. We have preservation of the norm of $\hat U_\ell$ if $\lambda$ is purely imaginary. Its norm decays if $\lambda$ has positive real part.

\subsection{Preservation of energy}

Consider now a situation in which all eigenvalues of $\hat A$ are purely imaginary. Then\footnote{The dagger denotes the Hermitian (conjugate) transpose.}
\begin{align}
 \hat U^\dagger \hat U = \hat U_{0}^\dagger \hat U_{0},
\end{align}
i.e., the $L^2$ norm of $\hat U$ is preserved. The $L^2$ norm of $\hat u$ is \emph{not} preserved, but it needs to be weighted:
\begin{align}
 \hat U^\dagger \hat U = \hat u^\dagger \hat R^\dagger \hat R \hat u.
\end{align}
Call $\hat M := \hat R^\dagger \hat R$. One observes both the symmetry $\hat M^\dagger = \hat M$ and the skew-symmetry of $\hat M \hat A$:
\begin{align}
\hat M \hat A = \hat R^\dagger \hat R \hat A \underbrace{\hat R^{-1} \hat R}_{=\text{id}} = \hat R^\dagger \hat \Lambda \hat R = -(\hat R^\dagger \hat \Lambda \hat R)^\dagger= -(\hat R^\dagger \hat R \hat A)^\dagger = -(\hat M \hat A)^\dagger,
\end{align}
where we have used $\Lambda^\dagger = -\Lambda$.
Thus, the Fourier symbol of the matrix $M$ is simply $\hat R^\dagger \hat R$, with as many free parameters as the dimension of the space on which $\hat M$ acts, one for each free scaling parameter of the eigenvector.

Consider now $A = D$, the space discretization of central Active Flux. $\hat D$ reads
\begin{align}
\left( \begin{array}{cc} \frac{1}{\tau} - \tau & 3 - \frac{3}{\tau} \\ \tau - 1 & 0 \end{array} \right)
\end{align}
and its eigenvalues are
\begin{align}
    \lambda^{\hat D}_\pm = \frac{\tau-1}{2\tau}\left(-\tau-1 \pm \sqrt{1 + 14 \tau + \tau^2}\right).
\end{align}
They are indeed imaginary, since upon the transformation $\tau \mapsto \frac{1}{\tau}$ (complex conjugation) they map into their negative
\begin{align}
    \frac{\frac{1}{\tau}-1}{2}\tau\left(-\frac{1}{\tau}-1 \pm \sqrt{1 + 14 \frac{1}{\tau} + \frac{1}{\tau^2}}\right) = -\frac{\tau-1}{2\tau}\left(-1-\tau \pm \sqrt{\tau^2 + 14 \tau + 1}\right) = - \lambda^{\hat D}_\pm.
\end{align}
An ingenious choice of normalization for the eigenvectors gives rise to the following Fourier symbol of the mass matrix:
\begin{align}
\hat M = \left( \begin{array}{cc} \alpha & \beta \\ \beta \tau & 3 \alpha + \beta (1+\tau) \end{array} \right ).
\end{align}
Here, $\alpha$ and $\beta$ are any Laurent polynomials in $\tau$ under the condition that the matrix remains Hermitian. This is equivalent to
\begin{align}
 \bar \beta &= \beta \tau, \label{eq:conditionbeta}\\
 \bar \alpha &= \alpha.
\end{align}
If $\beta = \sum_{\ell = -s}^s \beta_\ell \tau^\ell$, then condition \eqref{eq:conditionbeta} reads
\begin{align}
\bar \beta = \sum_{\ell = -s}^s \beta_\ell \tau^{-\ell} = \sum_{\ell = -s}^s \beta_{-\ell} \tau^{\ell} &= \sum_{\ell = -s+1}^{s+1} \beta_{\ell-1} \tau^{\ell} = \sum_{\ell = -s}^s \beta_\ell \tau^{\ell+1} =  \beta \tau.
\end{align}
The first polynomials for $s = 0, 1, 2, \dots$ that fulfill this are, up to scaling and linear combinations,
\begin{align}
0, 1 + \frac{1}{\tau}, \tau + \frac{1}{\tau^2}, \ldots
\end{align}

Similarly, $\alpha$ must be symmetric upon the transformation $\tau \mapsto \frac{1}{\tau}$, i.e. the first few polynomials (up to scaling and linear combinations) are
\begin{align}
1, \tau + \frac{1}{\tau}, \tau^2 + \frac{1}{\tau^2}, \dots
\end{align}

One choice is $\beta = 0, \alpha = 1$, yielding matrix \eqref{eq:central-AF-diagonal-mass}, another is
\begin{align}
 \alpha &=
 m_p + \frac{3 m_p - m_v + 2 m_{vv}}{6}\left( \tau + \frac{1}{\tau} \right)
 \in \mathrm{span}\left(1, \tau + \frac{1}{\tau}\right), \\
 \beta &= -(3 m_p - m_v)\frac{1+\tau}{2 \tau} \in \mathrm{span}\left(  1 + \frac{1}{\tau} \right).
\end{align}
yielding \eqref{eq:central-AF-banded-mass}.
Further choices can thus be systematically derived, for instance using
\begin{align*}
 \alpha &=
 m_p + y\left( \tau + \frac{1}{\tau} \right) + \frac{m_{vvv} - m_{vvp}}{3} \left( \tau^2 + \frac{1}{\tau^2} \right )
 \in \mathrm{span}\left(1, \tau + \frac{1}{\tau}, \tau^2 + \frac{1}{\tau^2}\right), \\
 \beta &= \frac{-3 m_p + m_v}{2}\left( 1 + \frac{1}{\tau}\right ) + m_{vvp} \left( \tau + \frac{1}{\tau^2} \right ) \in \mathrm{span}\left(  1 + \frac{1}{\tau}, \tau + \frac{1}{\tau^2}  \right).
\end{align*}
with $y := \frac{3 m_p - m_v + 2 m_{vv} - 2 m_{vvp}}{6}$
one obtains{\tiny
\begin{align*}
{\text{\normalsize $M =$}} \left( \begin{array}{c|cc|cc|cc|cc|cc|cc|c}
& &&& \text{\rotatebox{20}{$\ddots$}} &&&& \\\cline{1-11}
&\frac{m_{vvv} - m_{vvp}}{3} & m_{vvp} & y & \frac{m_v-3m_p}{2} & m_p & \frac{m_v-3m_p}{2} & y & m_{vvp} & \frac{m_{vvv} - m_{vvp}}{3} &  \\
& & m_{vvv} & m_{vvp} & m_{vv} & \frac{m_v-3m_p}{2} & m_v &\frac{m_v-3m_p}{2}  & m_{vv} & m_{vvp} & m_{vvv}
\\\cline{2-13}
\multicolumn{3}{l|}{}&\frac{m_{vvv} - m_{vvp}}{3} & m_{vvp} & y & \frac{m_v-3m_p}{2} & m_p & \frac{m_v-3m_p}{2} & y & m_{vvp} & \frac{m_{vvv} - m_{vvp}}{3} &  \\
\multicolumn{3}{l|}{}& & m_{vvv} & m_{vvp} & m_{vv} & \frac{m_v-3m_p}{2} & m_v &\frac{m_v-3m_p}{2}  & m_{vv} & m_{vvp} & m_{vvv} \\
\cline{4-14}
\multicolumn{5}{l|}{}&&&&& \text{\rotatebox{20}{$\ddots$}} &&&&
\end{array} \right)
\end{align*}}

\printbibliography

@misc{barsukow2025stabilityRepro,
  title={Reproducibility repository for
         "{S}tability of the Active Flux Method in the
         Framework of Summation-by-Parts Operators"},
  author={Barsukow, Wasilij and Klingenberg, Christian and
          Lechner, Lisa and Nordstr{\"o}m, Jan and Ortleb, Sigrun
          and Ranocha, Hendrik},
  year={2025},
  howpublished={\url{https://github.com/ranocha/2025_active_flux_sbp}},
  doi={10.5281/zenodo.15861045}
}

@article{svard2004coordinate,
  title={On Coordinate Transformations for Summation-by-Parts Operators},
  author={Sv{\"a}rd, Magnus},
  journal={Journal of Scientific Computing},
  volume={20},
  number={1},
  pages={29--42},
  year={2004},
  publisher={Springer},
  doi={10.1023/A:1025881528802}
}

@article{bezanson2017julia,
  title={Julia: {A} Fresh Approach to Numerical Computing},
  author={Bezanson, Jeff and Edelman, Alan and Karpinski, Stefan and
          Shah, Viral B},
  journal={SIAM Review},
  volume={59},
  number={1},
  pages={65--98},
  year={2017},
  publisher={SIAM},
  eprint={1411.1607},
  eprinttype={arxiv},
  eprintclass={cs.MS},
  doi={10.1137/141000671}
}

@article{danisch2021makie,
  title={Makie.jl: Flexible high-performance data visualization for {J}ulia},
  author={Danisch, Simon and Krumbiegel, Julius},
  journal={Journal of Open Source Software},
  volume={6},
  number={65},
  pages={3349},
  year={2021},
  doi={10.21105/joss.03349}
}

@article{abgrall2023combination,
  title={A combination of residual distribution and the active flux formulations
         or a new class of schemes that can combine several writings of the same
         hyperbolic problem: application to the {1D} {E}uler equations},
  author={Abgrall, R{\'e}mi},
  journal={Communications on Applied Mathematics and Computation},
  volume={5},
  number={1},
  pages={370--402},
  year={2023},
  publisher={Springer},
  doi={10.1007/s42967-021-00175-w}
}

@article{abgrall2023extensions,
  title={Extensions of active flux to arbitrary order of accuracy},
  author={Abgrall, Remi and Barsukow, Wasilij},
  journal={ESAIM: Mathematical Modelling and Numerical Analysis},
  volume={57},
  number={2},
  pages={991--1027},
  year={2023},
  publisher={EDP Sciences},
  doi={10.1051/m2an/2023004}
}

@article{abgrall2025semi,
  title={A Semi-discrete Active Flux Method for the {E}uler Equations on
         {C}artesian Grids},
  author={Abgrall, R{\'e}mi and Barsukow, Wasilij and Klingenberg, Christian},
  journal={Journal of Scientific Computing},
  volume={102},
  number={2},
  pages={36},
  year={2025},
  publisher={Springer},
  doi={10.1007/s10915-024-02749-1}
}

@article{fernandez2014review,
  title={Review of summation-by-parts operators with simultaneous approximation
         terms for the numerical solution of partial differential equations},
  author={Fern{\'a}ndez, David C Del Rey and Hicken, Jason E and Zingg, David W},
  journal={Computers {\&} Fluids},
  volume={95},
  pages={171--196},
  year={2014},
  publisher={Elsevier},
  doi={10.1016/j.compfluid.2014.02.016}
}

@article{svard2014review,
  title={Review of summation-by-parts schemes for initial-boundary-value problems},
  author={Sv{\"a}rd, Magnus and Nordstr{\"o}m, Jan},
  journal={Journal of Computational Physics},
  volume={268},
  pages={17--38},
  year={2014},
  publisher={Elsevier},
  doi={10.1016/j.jcp.2014.02.031}
}

@article{mattsson2017diagonal,
  title={Diagonal-norm upwind {SBP} operators},
  author={Mattsson, Ken},
  journal={Journal of Computational Physics},
  volume={335},
  pages={283--310},
  year={2017},
  publisher={Elsevier},
  doi={10.1016/j.jcp.2017.01.042}
}

@article{mattsson2007high,
  title={High-order accurate computations for unsteady aerodynamics},
  author={Mattsson, Ken and Sv{\"a}rd, Magnus and Carpenter, Mark and
          Nordstr{\"o}m, Jan},
  journal={Computers \& Fluids},
  volume={36},
  number={3},
  pages={636--649},
  year={2007},
  publisher={Elsevier},
  doi={10.1016/j.compfluid.2006.02.004}
}

@article{ortleb2023stability,
title={On the Stability of {IMEX} Upwind {gSBP} Schemes for {1D} Linear Advection-Diffusion Equations},
author={Ortleb, Sigrun},
journal={Communications on Applied Mathematics and Computation},
year={2023},
publisher={Springer},
doi={10.1007/s42967-023-00296-4}
}

@article{ranocha2021broadclass,
 title = {A Broad Class of Conservative Numerical Methods for Dispersive Wave Equations},
 author = {Hendrik Ranocha and Dimitrios Mitsotakis and David I. Ketcheson},
 journal = {Communications in Computational Physics},
 volume = {29},
 year = {2021},
 issue = {4},
 pages = {979--1029},
 doi = {10.4208/cicp.OA-2020-0119}
}

@article{stiernstroem2021aresidualbased,
title = {A residual-based artificial viscosity finite difference method for scalar conservation laws},
journal = {Journal of Computational Physics},
volume = {430},
pages = {110100},
year = {2021},
doi = {https://doi.org/10.1016/j.jcp.2020.110100},
author = {Vidar Stiernström and Lukas Lundgren and Murtazo Nazarov and Ken Mattsson},
}

@article{svard2005steady,
  title={Steady-state computations using summation-by-parts operators},
  author={Sv{\"a}rd, Magnus and Mattsson, Ken and Nordstr{\"o}m, Jan},
  journal={Journal of Scientific Computing},
  volume={24},
  pages={79--95},
  year={2005},
  publisher={Springer},
  doi={10.1007/s10915-004-4788-2}
}

@article{dovgilovich2015high,
  title={High-accuracy finite-difference schemes for solving elastodynamic problems in curvilinear coordinates within multiblock approach},
  author={Dovgilovich, Leonid and Sofronov, Ivan},
  journal={Applied Numerical Mathematics},
  volume={93},
  pages={176--194},
  year={2015},
  publisher={Elsevier},
  doi={10.1016/j.apnum.2014.06.005}
}

@article{barsukow24affourier,
  title={Analysis of the multi-dimensional semi-discrete {A}ctive {F}lux method using the Fourier transform},
  author={Barsukow, Wasilij and Kern, Janina and Klingenberg, Christian and Lechner, Lisa},
  journal={Communications on Applied Mathematics and Computation},
  pages={1--49},
  year={2025},
  publisher={Springer}
}

@article{roe21,
  title={Designing {CFD} methods for bandwidth—A physical approach},
  author={Roe, Philip},
  journal={Computers \& Fluids},
  volume={214},
  pages={104774},
  year={2021},
  publisher={Elsevier}
}

@article{carpenter1994timestable,
  title={Time-Stable Boundary Conditions for Finite-Difference Schemes Solving
         Hyperbolic Systems: {M}ethodology and Application to High-Order Compact
         Schemes},
  author={Carpenter, Mark H and Gottlieb, David and Abarbanel, Saul},
  journal={Journal of Computational Physics},
  volume={111},
  number={2},
  pages={220--236},
  year={1994},
  publisher={Elsevier},
  doi={10.1006/jcph.1994.1057}
}

@article{olsson1995_I,
 author = {P. Olsson},
 journal = {Mathematics of Computation},
 number = {211},
 pages = {1035--1065},
 publisher = {American Mathematical Society},
 title = {Summation by parts, projections, and stability. I},
 volume = {64},
 year = {1995},
 doi = {10.2307/2153512}
}

@article{olsson1995_II,
 author = {P. Olsson},
 journal = {Mathematics of Computation},
 number = {212},
 pages = {1473--1493},
 publisher = {American Mathematical Society},
 title = {Summation by parts, projections, and stability. II},
 volume = {64},
 year = {1995},
 doi = {10.2307/2153366}
}

@article{lundgren2020anefficient,
title = {An efficient finite difference method for the shallow water equations},
journal = {Journal of Computational Physics},
volume = {422},
pages = {109784},
year = {2020},
issn = {0021-9991},
doi = {10.1016/j.jcp.2020.109784},
author = {Lukas Lundgren and Ken Mattsson},
}

@ARTICLE{Nordström2006375,
	author = {Nordström, Jan},
	title = {Conservative finite difference formulations, variable coefficients, energy estimates and artificial dissipation},
	year = {2006},
	journal = {Journal of Scientific Computing},
	volume = {29},
	number = {3},
	pages = {375 – 404},
	doi = {10.1007/s10915-005-9013-4},
	type = {Article},
}

@ARTICLE{Fisher2013353,
	author = {Fisher, Travis C. and Carpenter, Mark H. and Nordström, Jan and Yamaleev, Nail K. and Swanson, Charles},
	title = {Discretely conservative finite-difference formulations for nonlinear conservation laws in split form: Theory and boundary conditions},
	year = {2013},
	journal = {Journal of Computational Physics},
	volume = {234},
	number = {1},
	pages = {353 – 375},
	doi = {10.1016/j.jcp.2012.09.026},
	type = {Article},
}

@article{ortleb2017akinetic,
author = {Ortleb, Sigrun},
title = {A kinetic energy preserving {DG} scheme based on {G}auss-{L}egendre points},
journal = {Journal of Scientifc Computing},
year = {2017},
volume = {71},
pages = {1135--1168},
doi = {10.1007/s10915-016-0334-2}
}

@article{glaubitz2023summation,
author = {Glaubitz, Jan and Nordstr\"{o}m, Jan and \"{O}ffner, Philipp},
title = {Summation-by-Parts Operators for General Function Spaces},
journal = {SIAM Journal on Numerical Analysis},
volume = {61},
number = {2},
pages = {733-754},
year = {2023},
doi = {10.1137/22M1470141}
}

@ARTICLE{Glaubitz2023_2,
	author = {Glaubitz, Jan and Klein, Simon-Christian and Nordström, Jan and Öffner, Philipp},
	title = {Multi-dimensional summation-by-parts operators for general function spaces: Theory and construction},
	year = {2023},
	journal = {Journal of Computational Physics},
	volume = {491},
	doi = {10.1016/j.jcp.2023.112370},
	type = {Article},
}

@ARTICLE{Glaubitz2024_1,
	author = {Glaubitz, Jan and Nordström, Jan and Öffner, Philipp},
	title = {Energy-Stable Global Radial Basis Function Methods on Summation-By-Parts Form},
	year = {2024},
	journal = {Journal of Scientific Computing},
	volume = {98},
	number = {1},
	doi = {10.1007/s10915-023-02427-8},
	type = {Article},
}

@ARTICLE{Glaubitz2024_2,
	author = {Glaubitz, Jan and Klein, Simon-Christian and Nordström, Jan and Öffner, Philipp},
	title = {Summation-by-parts operators for general function spaces: The second derivative},
	year = {2024},
	journal = {Journal of Computational Physics},
	volume = {504},
	doi = {10.1016/j.jcp.2024.112889},
	type = {Article},
}

@article{ranocha2025ontherobustness,
title = {On the robustness of high-order upwind summation-by-parts methods for nonlinear conservation laws},
journal = {Journal of Computational Physics},
volume = {520},
pages = {113471},
year = {2025},
issn = {0021-9991},
doi = {https://doi.org/10.1016/j.jcp.2024.113471},
author = {Hendrik Ranocha and Andrew R. Winters and Michael Schlottke-Lakemper and Philipp \"{O}ffner and Jan Glaubitz and Gregor J. Gassner}
}

@inproceedings{kreiss1974finite,
  title={Finite Element and Finite Difference Methods for Hyperbolic
         Partial Differential Equations},
  author={Kreiss, Heinz-Otto and Scherer, Godela},
  booktitle={Mathematical Aspects of Finite Elements in Partial
             Differential Equations},
  editor={de Boor, Carl},
  pages={195--212},
  year={1974},
  publisher={Academic Press},
  address={New York}
}

@article{strand1994summation,
  title={Summation by Parts for Finite Difference Approximations for
         {$d/dx$}},
  author={Strand, Bo},
  journal={Journal of Computational Physics},
  volume={110},
  number={1},
  pages={47--67},
  year={1994},
  publisher={Elsevier},
  doi={10.1006/jcph.1994.1005}
}

@article{nordstrom2001finite,
  title={Finite volume approximations and strict stability for hyperbolic problems},
  author={Nordstr{\"o}m, Jan and Bj{\"o}rck, Martin},
  journal={Applied Numerical Mathematics},
  volume={38},
  number={3},
  pages={237--255},
  year={2001},
  publisher={Elsevier},
  doi={10.1016/S0168-9274(01)00027-7}
}

@article{nordstrom2003finite,
  title={Finite volume methods, unstructured meshes and strict stability
         for hyperbolic problems},
  author={Nordstr{\"o}m, Jan and Forsberg, Karl and Adamsson, Carl and Eliasson, Peter},
  journal={Applied Numerical Mathematics},
  volume={45},
  number={4},
  pages={453--473},
  year={2003},
  publisher={Elsevier},
  doi={10.1016/S0168-9274(02)00239-8}
}

@article{hicken2016multidimensional,
  title={Multidimensional Summation-By-Parts Operators:
         {G}eneral Theory and Application to Simplex Elements},
  author={Hicken, Jason E and Fern{\'a}ndez, David C Del Rey and Zingg, David W},
  journal={SIAM Journal on Scientific Computing},
  volume={38},
  number={4},
  pages={A1935--A1958},
  year={2016},
  publisher={Society for Industrial and Applied Mathematics},
  doi={10.1137/15M1038360}
}

@article{hicken2020entropy,
  title={Entropy-stable, high-order summation-by-parts discretizations
         without interface penalties},
  author={Hicken, Jason E},
  journal={Journal of Scientific Computing},
  volume={82},
  number={2},
  pages={50},
  year={2020},
  publisher={Springer},
  doi={10.1007/s10915-020-01154-8}
}

@article{abgrall2020analysisI,
  title={Analysis of the {SBP-SAT} Stabilization for Finite Element Methods
         Part {I}: {L}inear problems},
  author={Abgrall, R{\'e}mi and Nordstr{\"o}m, Jan and {\"O}ffner, Philipp
          and Tokareva, Svetlana},
  journal={Journal of Scientific Computing},
  volume={85},
  number={2},
  pages={1--29},
  year={2020},
  publisher={Springer},
  doi={10.1007/s10915-020-01349-z},
  eprint={1912.08108},
  eprinttype={arxiv},
  eprintclass={math.NA}
}

@article{gassner2013skew,
  title={A Skew-Symmetric Discontinuous {G}alerkin Spectral Element
         Discretization and Its Relation to {SBP}-{SAT} Finite Difference
         Methods},
  author={Gassner, Gregor Josef},
  journal={SIAM Journal on Scientific Computing},
  volume={35},
  number={3},
  pages={A1233--A1253},
  year={2013},
  publisher={Society for Industrial and Applied Mathematics},
  doi={10.1137/120890144}
}

@article{carpenter2014entropy,
  title={Entropy Stable Spectral Collocation Schemes for the
         {N}avier-{S}tokes Equations: {D}iscontinuous Interfaces},
  author={Carpenter, Mark H and Fisher, Travis C and Nielsen, Eric J and
          Frankel, Steven H},
  journal={SIAM Journal on Scientific Computing},
  volume={36},
  number={5},
  pages={B835--B867},
  year={2014},
  publisher={Society for Industrial and Applied Mathematics},
  doi={10.1137/130932193}
}

@article{chan2018discretely,
  title={On discretely entropy conservative and entropy stable discontinuous
         {G}alerkin methods},
  author={Chan, Jesse},
  journal={Journal of Computational Physics},
  volume={362},
  pages={346--374},
  year={2018},
  publisher={Elsevier},
  doi={10.1016/j.jcp.2018.02.033}
}

@inproceedings{huynh2007flux,
  title={A Flux Reconstruction Approach to High-Order Schemes Including
         Discontinuous {G}alerkin Methods},
  author={Huynh, H. T.},
  booktitle={18th AIAA Computational Fluid Dynamics Conference},
  year={2007},
  organization={American Institute of Aeronautics and Astronautics},
  doi={10.2514/6.2007-4079}
}

@article{vincent2011newclass,
  title={A New Class of High-Order Energy Stable Flux Reconstruction Schemes},
  author={Vincent, Peter E and Castonguay, Patrice and Jameson, Antony},
  journal={Journal of Scientific Computing},
  volume={47},
  number={1},
  pages={50--72},
  year={2011},
  publisher={Springer},
  doi={10.1007/s10915-010-9420-z}
}

@article{ranocha2016summation,
  title={Summation-by-parts operators for correction procedure via
         reconstruction},
  author={Ranocha, Hendrik and {\"O}ffner, Philipp and Sonar, Thomas},
  journal={Journal of Computational Physics},
  volume={311},
  pages={299--328},
  year={2016},
  month={04},
  publisher={Elsevier},
  doi={10.1016/j.jcp.2016.02.009},
  eprint={1511.02052},
  eprinttype={arxiv},
  eprintclass={math.NA}
}

@online{hicken2024constructing,
  title={Constructing stable, high-order finite-difference operators
         on point clouds over complex geometries},
  author={Hicken, Jason and Yan, Ge and Kaur, Sharanjeet},
  year={2024},
  month={09},
  eprint={2409.00809},
  eprinttype={arxiv},
  eprintclass={math.NA}
}

@article{svard2019convergence,
  title={On the convergence rates of energy-stable finite-difference schemes},
  author={Sv{\"a}rd, Magnus and Nordstr{\"o}m, Jan},
  journal={Journal of Computational Physics},
  volume={397},
  pages={108819},
  year={2019},
  publisher={Elsevier},
  doi={10.1016/j.jcp.2019.07.018}
}

@article{trapp1973blockcirculant,
  title={Inverses of circulant matrices and block circulant matrices},
  author={Trapp, George E},
  journal={Kyungpook Mathematical Journal},
  year={1973},
  volume={13},
  number={1},
  pages={11--20}
}

@book{davis1979,
  title     = "Circulant Matrices",
  author    = "Davis, Philipp J.",
  year      = 1979,
  publisher = "Wiley Interscience",
  address   = "New York"
}

@article{ranocha2021optimized,
  title={Optimized {R}unge-{K}utta Methods with Automatic Step Size Control
         for Compressible Computational Fluid Dynamics},
  author={Ranocha, Hendrik and Dalcin, Lisandro and Parsani, Matteo
          and Ketcheson, David I.},
  journal={Communications on Applied Mathematics and Computation},
  volume={4},
  pages={1191--1228},
  year={2021},
  month={11},
  doi={10.1007/s42967-021-00159-w},
  eprint={2104.06836},
  eprinttype={arxiv},
  eprintclass={math.NA}
}

@article{ketcheson2019relaxation,
  title={Relaxation {R}unge-{K}utta Methods: {C}onservation and
         Stability for Inner-Product Norms},
  author={Ketcheson, David I},
  journal={SIAM Journal on Numerical Analysis},
  volume={57},
  number={6},
  pages={2850--2870},
  year={2019},
  publisher={Society for Industrial and Applied Mathematics},
  doi={10.1137/19M1263662},
  eprint={1905.09847},
  eprinttype={arxiv},
  eprintclass={math.NA}
}

@article{ranocha2020relaxation,
  title={Relaxation {R}unge-{K}utta Methods: Fully-Discrete Explicit
         Entropy-Stable Schemes for the Compressible {E}uler and
         {N}avier-{S}tokes Equations},
  author={Ranocha, Hendrik and Sayyari, Mohammed and Dalcin, Lisandro and
          Parsani, Matteo and Ketcheson, David I.},
  journal={SIAM Journal on Scientific Computing},
  volume={42},
  number={2},
  pages={A612--A638},
  year={2020},
  month={03},
  publisher={Society for Industrial and Applied Mathematics},
  doi={10.1137/19M1263480},
  eprint={1905.09129},
  eprinttype={arxiv},
  eprintclass={math.NA}
}

@article{ranocha2020general,
  title={General Relaxation Methods for Initial-Value Problems
         with Application to Multistep Schemes},
  author={Ranocha, Hendrik and L{\'o}czi, Lajos and Ketcheson, David I},
  journal={Numerische Mathematik},
  year={2020},
  month={10},
  volume={146},
  pages={875--906},
  publisher={Springer Nature},
  doi={10.1007/s00211-020-01158-4},
  eprint={2003.03012},
  eprinttype={arxiv},
  eprintclass={math.NA}
}

@article{eymann13,
author = {Timothy A. Eymann and Philip L. Roe},
title = {Multidimensional Active Flux Schemes},
journal = {21st AIAA Computational Fluid Dynamics Conference},
chapter = {},
year = {2013},
pages = {},
doi = {10.2514/6.2013-2940},
}

@article{eymann11,
	author = {Timothy A. Eymann and Philip L. Roe},
	title = {{Active Flux Schemes}},
	journal = {49th AIAA Aerospace Sciences Meeting including the New Horizons Forum and Aerospace Exposition},
	chapter = {},
	year = {2011},
	pages = {},
	doi = {10.2514/6.2011-382},
}

@article{vanleer77,
author = {Bram Van Leer},
title = {Towards the ultimate conservative difference scheme. IV. A new approach to numerical convection},
journal = {Journal of Computational Physics},
volume = {23},
number = {3},
pages = {276-299},
year = {1977},
issn = {0021-9991},
doi = {10.1016/0021-9991(77)90095-X},
}

@article{barsukow18activeflux,
author = {Wasilij Barsukow},
title = {Stationarity preserving schemes for multi-dimensional linear systems},
journal = {Mathematics of Computation},
volume = {88},
number = {},
pages = {1621-1645},
year = {2019},
doi={10.1090/mcom/3394}
}

@article{abgrall23,
	author = {{Abgrall, Remi} and {Barsukow, Wasilij}},
	title = {Extensions of Active Flux to arbitrary order of accuracy},
	DOI= "10.1051/m2an/2023004",
	journal = {ESAIM: M2AN},
	year = 2023,
	volume = 57,
	number = 2,
	pages = "991-1027",
}

@article{abgrall23a,
title = {A hybrid finite element - finite volume method for conservation laws},
journal = {Applied Mathematics and Computation},
volume = {447},
pages = {},
year = {2023},
issn = {0096-3003},
doi = {https://doi.org/10.1016/j.amc.2023.127846},

author = {R√©mi Abgrall and Wasilij Barsukow}
}

@article{Lechner25,
  title={A generalized Active Flux method of arbitrarily high order in two dimensions},
  author={Barsukow, Wasilij and Chandrashekar, Praveen and Klingenberg, Christian and Lechner, Lisa},
  journal={arXiv preprint arXiv:2502.05101},
  year={2025}
}

@article{eymann2011active,
  title={Active flux schemes for systems},
  author={Eymann, Timothy and Roe, Philip},
  journal={20th AIAA computational fluid dynamics conference},
   volume = {AIAA 2011-3840},
  pages={},
  year={2011},
  doi={10.2514/6.2011-3840}
}

@book{hairer2006geometric,
  title={Geometric Numerical Integration: {S}tructure-Preserving
         Algorithms for Ordinary Differential Equations},
  author={Hairer, Ernst and Lubich, Christian and Wanner, Gerhard},
  series={Springer Series in Computational Mathematics},
  volume={31},
  year={2006},
  publisher={Springer-Verlag},
  address={Berlin Heidelberg},
  doi={10.1007/3-540-30666-8}
}

@article{sun2017stability,
  title={Stability of the fourth order {R}unge-{K}utta method for
         time-dependent partial differential equations},
  author={Sun, Zheng and Shu, Chi-Wang},
  journal={Annals of Mathematical Sciences and Applications},
  volume={2},
  number={2},
  pages={255--284},
  year={2017},
  doi={10.4310/AMSA.2017.v2.n2.a3}
}

@article{sun2019strong,
  title={Strong Stability of Explicit {R}unge-{K}utta Time
         Discretizations},
  author={Sun, Zheng and Shu, Chi-Wang},
  journal={SIAM Journal on Numerical Analysis},
  volume={57},
  number={3},
  pages={1158--1182},
  year={2019},
  publisher={SIAM},
  doi={10.1137/18M122892X},
  eprint={1811.10680},
  eprinttype={arxiv},
  eprintclass={math.NA}
}

@incollection{tadmor2002semidiscrete,
  title={From Semidiscrete to Fully Discrete: {S}tability of
         {R}unge-{K}utta Schemes by the Energy Method {II}},
  author={Tadmor, Eitan},
  booktitle={Collected Lectures on the Preservation of Stability under Discretization},
  editor={Estep, Donald J and Tavener, Simon},
  series={Proceedings in Applied Mathematics},
  volume={109},
  pages={25--49},
  year={2002},
  publisher={Society for Industrial and Applied Mathematics},
  address={Philadelphia}
}

@article{ranocha2018L2stability,
  title={{$L_2$} Stability of Explicit {R}unge-{K}utta Schemes},
  author={Ranocha, Hendrik and {\"O}ffner, Philipp},
  journal={Journal of Scientific Computing},
  volume={75},
  number={2},
  pages={1040--1056},
  year={2018},
  month={05},
  doi={10.1007/s10915-017-0595-4}
}

@article{ranocha2021strong,
  title={On Strong Stability of Explicit {R}unge-{K}utta Methods for
         Nonlinear Semibounded Operators},
  author={Ranocha, Hendrik},
  journal={IMA Journal of Numerical Analysis},
  year={2021},
  month={01},
  volume={41},
  number={1},
  pages={654--682},
  publisher={Oxford University Press},
  doi={10.1093/imanum/drz070},
  eprint={1811.11601},
  eprinttype={arxiv},
  eprintclass={math.NA}
}

@article{ranocha2020energy,
  title={Energy Stability of Explicit {R}unge-{K}utta Methods for
         Nonautonomous or Nonlinear Problems},
  author={Ranocha, Hendrik and Ketcheson, David I},
  journal={SIAM Journal on Numerical Analysis},
  year={2020},
  month={11},
  volume={58},
  number={6},
  pages={3382--3405},
  publisher={Society for Industrial and Applied Mathematics},
  doi={10.1137/19M1290346},
  eprint={1909.13215},
  eprinttype={arxiv},
  eprintclass={math.NA}
}

@article{lozano2018entropy,
  title={Entropy Production by Explicit {R}unge-{K}utta Schemes},
  author={Lozano, Carlos},
  journal={Journal of Scientific Computing},
  volume={76},
  number={1},
  pages={521--565},
  year={2018},
  publisher={Springer},
  doi={10.1007/s10915-017-0627-0}
}

@article{mattsson2013solution,
  title={A solution to the stability issues with block norm summation by parts operators},
  author={Mattsson, Ken and Almquist, Martin},
  journal={Journal of Computational Physics},
  volume={253},
  pages={418--442},
  year={2013},
  publisher={Elsevier},
  doi={10.1016/j.jcp.2013.07.013}
}

\end{document}